%% file: main-journal.tex
\documentclass[a4paper,10pt]{article}

\usepackage{amsfonts,amssymb,amsmath}
\usepackage{graphicx,graphics,epsfig,color}
\usepackage{multicol}
\usepackage{rgsMacros}
\usepackage{graphicx} 
\usepackage{array,amssymb,color,xspace,psfrag}
\usepackage[numbers,sort&compress,square]{natbib} 
\usepackage{amssymb,enumerate}
\usepackage{flushend}

\usepackage{hyperref}

\usepackage{enumitem}
\usepackage{balance}
\usepackage{wrapfig}

    \definecolor{gray}{rgb}{0.33,0.4,0.47}
    \definecolor{steelblue}{rgb}{0,.42,.7}
    \definecolor{britishgreen}{rgb}{0,0.26,0.15}
    \definecolor{navyblue}{rgb}{0,0,.8}
    \definecolor{olivegreen}{rgb}{0.14,0.29,0}
    \definecolor{myred}{rgb}{0.86,0.1,0.16}

 \newif\ifitsdraft
 \def\itsdraft{\global\itsdrafttrue}

 \itsdraft  

\ifitsdraft

\usepackage{fullpage}
\usepackage{multirow}
\usepackage{subfig}
\usepackage{hyperref}
\usepackage{booktabs}
\usepackage{xspace}
\usepackage[numbers,sort&compress]{natbib}
\usepackage{pstool}
\usepackage{psfrag}
\newtheorem{definition}{Definition}
\newtheorem{lemma}{Lemma}
\newtheorem{claim}{Claim}
\newtheorem{theorem}{Theorem}
\newtheorem{corollary}{Corollary}
\newtheorem{proposition}{Proposition}
\newtheorem{assumption}{Assumption} 

\newtheorem{remark}{Remark}
\newtheorem{example}{Example}
\newenvironment{proof}{\noindent {\it Proof.}}{ \hfill \mbox{\footnotesize $\blacksquare$ } \\ }

\allowdisplaybreaks
\pdfoutput = 1

\title{Sufficient conditions for forward invariance and contractivity in hybrid inclusions using barrier functions}    
\date{}

\author{Mohamed Maghenem and Ricardo G. Sanfelice
\thanks{ \blue{ M. Maghenem and R. G. Sanfelice are with Department of Electrical and Computer Engineering, University of California, Santa Cruz. Email: mmaghene,ricardo@ucsc.edu. 
This research has been partially supported by the National Science Foundation under Grant no. ECS-1710621 and Grant no. CNS-1544396, by the Air Force Office of Scientific Research under Grant no. FA9550-16-1-0015, Grant no. FA9550-19-1-0053, and Grant no. FA9550-19-1-0169, and by CITRIS and the Banatao Institute at the University of California.} } }

\begin{document}
\maketitle 
\begin{abstract} 
This paper studies set invariance and contractivity in hybrid systems modeled by hybrid inclusions using barrier functions. After introducing the notion of barrier function, we investigate sufficient conditions to guarantee different forward invariance and contractivity notions of a closed set for hybrid systems with nonuniqueness of solutions and solutions terminating prematurely. More precisely, we consider forward (pre-)invariance of sets, which guarantees that the maximal solutions starting from the set stay in it, and (pre-)contractivity, which further requires that the solutions starting from the boundary of the set evolve immediately (continuously or discretely) towards its interior. Our conditions for forward invariance and contractivity are infinitesimal and in terms of the proposed barrier functions. Examples illustrate the results.
\end{abstract}

\else

\usepackage{lipsum}

\newcommand\blfootnote[1]{%
  \begingroup
  \renewcommand\thefootnote{}\footnote{#1}%
  \addtocounter{footnote}{-1}%
  \endgroup}

\newtheorem{definition}{Definition}
\newtheorem{lemma}{Lemma}
\newtheorem{theorem}{Theorem}
\theoremstyle{plain}
\newtheorem{corollary}{Corollary}
\newtheorem{proposition}{Proposition}
\newtheorem{assumption}{Assumption} 
\newtheorem{remark}{Remark}
\newtheorem{condition}{Condition}
\newtheorem{example}{Example}
\begin{document}
\begin{frontmatter}
\title{Sufficient conditions for forward invariance and contractivity in hybrid inclusions using barrier functions}
\author[ucsc]{Mohamed Maghenem} and \
\ead{mmaghene@ucsc.edu}
\author[ucsc]{Ricardo G. Sanfelice} \
\ead{ricardo@ucsc.edu}
\address[ucsc]{Department of Electrical and Computer Engineering, University of California, Santa Cruz.}
\maketitle


\begin{abstract} 
This paper studies set invariance and contractivity in hybrid systems modeled by hybrid inclusions using barrier functions. After introducing the notion of barrier function, we investigate sufficient conditions to guarantee different forward invariance and contractivity notions of a closed set for hybrid systems with nonuniqueness of solutions and solutions terminating prematurely. We consider forward (pre-)invariance of sets, which guarantees that the maximal solutions starting from the set stay in it, and (pre-)contractivity, which further requires that the solutions starting from the boundary of the set evolve immediately (continuously or discretely) towards its interior. Our conditions for forward invariance and contractivity are infinitesimal and in terms of the proposed barrier functions. Examples illustrate the results. 
\end{abstract}
\end{frontmatter}
\fi

\section{Introduction}
Forward invariance of sets for dynamical systems is a property that requires the solutions starting in the considered set to remain in it along their entire domain of definition. The main challenge when studying forward invariance consists of providing sufficient conditions while avoiding explicit computation of the system solutions.\ifitsdraft\else\blfootnote{This research has been partially supported by the National Science Foundation under Grant no. ECS-1710621 and Grant no. CNS-1544396, by the Air Force Office of Scientific Research under Grant no. FA9550-16-1-0015, Grant no. FA9550-19-1-0053, and Grant no. FA9550-19-1-0169, and by CITRIS and the Banatao Institute at the University of California.} \fi

The study of set invariance for dynamical systems is a key step towards analyzing their stability and safety properties. Indeed, forward invariance has a close relationship to safety, which is a property that requires the system solutions starting from a given set of initial conditions to remain in a desired safe region \cite{wieland2007constructive}. Safety, also named conditional invariance in \cite{ladde1974}, is equivalent to forward invariance of a set, known as \textit{inductive invariant} in \cite{taly2009deductive}, which contains the set of initial conditions and not intersecting with the unsafe region \cite{prajna2005optimization}. Furthermore, the study of set invariance can be a key step to conclude stability properties for the system via relaxed Lyapunov conditions, including the well-known \textit{invariance principle} \cite{sanfelice2007invariance} and \textit{Matrosov Theorem} \cite{sanfelice2009asymptotic}. In addition to safety, the study of set invariance has also been extended in order to guarantee some closely related notions such as \textit{quasi-invariance}, \textit{conditional quasi-invariance} \cite{ladde1972analysis}, and \textit{contractivity} \cite{Blanchini:1999:SPS:2235754.2236030}. The contractivity property is a strong form of forward invariance. Indeed, a contractive set is forward invariant, and whenever a solution starts from its boundary, it immediately leaves the boundary and evolves towards its interior. Contractivity is also named \textit{strict invariance} in \cite{Aubin:1991:VT:120830}. Analyzing contractivity, which is a strong form of forward invariance, in the context of barrier functions is useful. Indeed, contractivity is closely related to invariance and it has been studied in the literature, for example, when computing \textit{set-induced Lyapunov functions} \cite{fiacchini2015computation, 272351, giesl2015review}. Hence, our goal is to formalize the differences between the two notions and how to certify each of them using barrier functions. 
 
\subsection{Background}

The interest in the study and characterization of forward invariance, while avoiding the computation of the system's solutions, dates back to the seminal work of Nagumo in \cite{nagumo1942lage}. In this reference, conditions involving the contingent cone and the system's dynamics on the boundary of a closed set are shown to be necessary and sufficient to conclude, from each point in the set, the existence of at least one solution that remains in the set. This last property defines what is known as \textit{weak forward invariance} \cite{clarke2008nonsmooth}, which is equivalent to forward invariance when the system's solutions are unique. Extensions of this result, using a similar type of cone conditions, are presented in \cite{Aubin:1991:VT:120830} to conclude weak forward invariance\footnote{Weak forward invariance is named \textit{viability} in \cite{Aubin:1991:VT:120830}.} for differential inclusions, in \cite{aubin2002impulse} for impulse differential inclusions, and in \cite{chai2015notions} for hybrid inclusions. For systems with continuous-time dynamics, when all the solutions starting from a closed set are required to remain in it, as stressed in \cite{Aubin:1991:VT:120830}, the invariance conditions concern the system's dynamics outside the set rather than on its boundary. As a consequence, the external contingent cone is introduced and used in \cite{Aubin:1991:VT:120830}. 

The relatively stronger form of invariance called contractivity is characterized in \cite{Blanchini:1999:SPS:2235754.2236030} in terms of the \textit{Minkowski} functional, both for differential and difference equations, for the particular case of convex and compact sets. For general closed sets and for systems modeled as differential inclusions, sufficient conditions are proposed in \cite{Aubin:1991:VT:120830} using, roughly speaking, the interior of the contingent cone and the system's dynamics at the boundary of the considered closed set. In general, the computation of the tangent cones is not a trivial task. However, when the considered set is defined as the zero-sublevel of a vector function, named \textit{barrier function candidate}, it is possible, under appropriate assumptions, to formulate invariance and contractivity conditions using only the barrier function candidate and the system's dynamics. 
 
The latter approach has been adopted in \cite{ames2014control} and \cite{glotfelter2017nonsmooth, 8625554} for differential equations and inclusions, respectively, and in \cite{prajna2007framework} for hybrid automata. In \cite{ames2014control}, another type of barrier function candidate is also considered. Such a barrier function candidate is positive and locally bounded on the interior of the considered set, and approaches infinity as its argument converges to the boundary of the set. When using the latter notion of barrier functions, solutions starting from the interior of the set to render invariant are not allowed to reach its boundary. In \cite{glotfelter2017nonsmooth}, sufficient conditions for invariance in terms of nonsmooth but locally Lipschitz barrier functions are considered with application to safe navigation for networks of vehicles in the presence of obstacles. See also \cite{xu2018correctness, nguyen2015safety} for more applications. Finally, in \cite{prajna2007framework}, methods to synthesize barrier functions are investigated. The latter work is extended in \cite{DAI201762, 10.1007/978-3-642-39799-8_17} by relaxing the conditions constraining the continuous-time evolution of the hybrid automata.

\subsection{Motivation}

A hybrid inclusion is defined as a differential inclusion with a constraint, which models the {\em flow} or continuous evolution of the system, and a difference inclusion with a constraint, modeling the {\em jumps}, or discrete events. In particular, handling nonuniqueness of solutions in hybrid inclusions and solutions terminating prematurely lead to particular forms of forward invariance and contractivity properties that we call forward pre-invariance and pre-contractivity, respectively, where the prefix ``pre'' indicates that some solutions may have a bounded (hybrid) time domain. The aforementioned notions have not been covered in the literature using barrier functions. Furthermore, having sufficient conditions for forward invariance in terms of barrier functions is useful, especially when control inputs are used to force such conditions \cite{wieland2007constructive, ames2014control}, or when the invariance, or the contractivity, task is to be combined with a stabilization task to be achieved inside a safety set \cite{romdlony2016stabilization}. Furthermore, in many applications, it is often the case that the closed set to be rendered forward invariant or contractive corresponds to the region where multiple scalar functions are nonpositive simultaneously. In such a case, it is typically difficult to find a single scalar function that defines the set of interest and, at the same time, is sufficiently smooth. This fact motivates the development of sufficient conditions guaranteeing forward invariance and contractivity when multiple scalar candidates define the considered set.
 
\subsection{Contributions}

In this paper, we introduce barrier functions and tools to certify forward invariance and contractivity in hybrid systems modeled as hybrid inclusions. We define a barrier function candidate as a vector function of the state variables. Sufficient conditions in terms of infinitesimal inequalities -- namely, without using information about solutions -- are proposed to guarantee that the set of points on which all the components of the barrier function candidate are nonpositive is forward invariant or contractive. More precisely, under mild conditions on the data defining the hybrid inclusion, we present conditions such that a barrier function candidate guarantees forward pre-invariance. The proposed conditions can be decomposed into flow and jump conditions that restrict the continuous and the discrete evolution of the hybrid system, respectively.

In Section \ref{Sec.2}, sufficient conditions for forward invariance that apply when the barrier function candidate is continuously differentiable and only locally Lipschitz are presented; see, respectively, Theorems \ref{prop1} and \ref{prop6bis}. 
Furthermore, the flow conditions therein need to be satisfied on a neighborhood outside the considered set. The flow conditions in Theorems \ref{prop1} and \ref{prop6bis} are relaxed in Proposition \ref{pro.nuniq} and Remark \ref{remns} using uniqueness functions. This relaxation is original to the best of our knowledge.  On the other hand, under a regularity condition on the gradient of the barrier function candidate, known as \textit{transversality condition}, which is typically assumed in the literature (see, e.g.,\cite{clarke1990optimization, Aubin:1991:VT:120830}), plus some extra regularity conditions on the flow dynamics, in Theorem \ref{propnou}, the flow condition needs to hold only on a smaller region of the state space which is the boundary of the considered set.

In Section \ref{Sec.4}, we analyze contractivity properties for hybrid systems. After a brief overview of the case of sets that are compact and convex, which itself extends to hybrid inclusions the results in \cite{Blanchini:1999:SPS:2235754.2236030} for differential and difference equations, we introduce a notion of contractivity for general closed sets. Furthermore, it extends to hybrid inclusions what is proposed in \cite{Aubin:1991:VT:120830} for differential inclusions only. The proposed notion essentially requires the system's solutions to evolve from points on the boundary of $K$ towards its interior via a flow or a jump. Sufficient conditions for contractivity in terms of the barrier function candidates defining the set are established when the latter candidate is either continuously differentiable or only locally Lipschitz, see Theorems \ref{prop5} and \ref{prop5lip}.

The results in this paper extend what was proposed in \cite{ames2014control, prajna2007framework, DAI201762, 10.1007/978-3-642-39799-8_17, glotfelter2017nonsmooth} to the more general context of hybrid systems modeled by hybrid inclusions. That is, hybrid inclusions offer many technical challenges that have not been handled in the existing literature. Those challenges are mainly due to the fact that the continuous-time evolution of the hybrid inclusion is not necessarily defined on an open set. Moreover, the considered set is defined as the zero-sublevel set of a barrier function candidate restricted to the set where the dynamics are defined; namely, the union of the flow and the jump sets. Since the latter sets can be closed, elements around the intersection between the two boundaries; namely, the zero-level set of the barrier candidate and the boundary of the flow set, needs a particular treatment. Our sufficient conditions using barrier functions are alternatives to those proposed in \cite{chai2015notions} and \cite{Aubin:1991:VT:120830} using tangent cone-based conditions. Indeed, our conditions exploit the fact that the set is the intersection of zero-sublevel sets of scalar functions; hence, the obtained conditions avoid as much as possible the computation of tangent cones. The latter task is known to be numerically expensive in some cases. It is also to be noted that some of our results build upon the well-known cone-based conditions in \cite{Aubin:1991:VT:120830} and    \cite{clarke2008nonsmooth}.

To the best of our knowledge, this is the first time in the literature where the concept of barrier functions is used for hybrid inclusions to analyze different set-invariance properties. Preliminary versions of this work are in the conference papers \cite{CP5-SIMUL3-CDC2018, CP5-SIMUL3-ACC2019}. However, only scalar barrier functions are considered in \cite{magh2018barrier}, and many proofs, explanations, and examples are omitted in both submissions. Furthermore, compared with \cite{CP5-SIMUL3-CDC2018, CP5-SIMUL3-ACC2019}, new results are proposed in Theorem \ref{propnou}, Theorem \ref{prop6}, Theorem \ref{prop6bis}, Proposition \ref{pro.nuniq}, and Theorem \ref{prop5lip}. 

The remainder of the paper is organized as follows. Preliminaries and basic conditions are presented in Section \ref{Sec.1}. Sufficient conditions of forward pre-invariance and invariance using barrier functions are in Section \ref{Sec.2}. Sufficient conditions for pre-contractivity and contractivity using barrier functions are in Section \ref{Sec.4}, respectively. Examples are included at each step in order to illustrate the proposed statements.
\ifitsdraft \else More technical details can be found in \cite{draftautomatica}.\fi

\textbf{Notation.} Let $\mathbb{R}_{\geq 0} := [0, \infty)$, 
$\mathbb{N} := \left\{0,1,\ldots \right\}$, and 
$\mathbb{N}^*:= \left\{1, 2,  \ldots, \infty \right\}$. For $x$, $y \in \mathbb{R}^n$ and a nonempty set $K \subset \mathbb{R}^n$, $x^\top$ denotes the transpose of $x$, $|x|$ the norm of $x$, $|x|_K := \inf_{y \in K} |x-y|$ defines the distance between $x$ and the set $K$,  $\langle x,y \rangle = x^\top y$ denotes the inner product between $x$ and $y$, and $\langle x,K \rangle = x^\top K := \defset{ x^\top z}{z \in K}$. The inequalities $x \leq 0$ and $x < 0$ mean that $x_i \leq 0$ and, respectively, $x_i < 0$ for all $i \in \left\{ 1,2,\ldots,n \right\}$. The opposites, namely, $x \nleq 0$ and $x \nless 0$ mean that there exists $i \in \left\{ 1,2,\ldots,n \right\}$ such that  $x_i > 0$ and, respectively, $x_i \geq 0$. For a set $K \subset \mathbb{R}^n$, we use $\mbox{int} (K)$ to denote its interior, $\partial K$ to denote its boundary, $\mbox{cl}(K)$ to denote its closure, and $U(K)$ to denote an open neighborhood around $K$. For a set $O \subset \mathbb{R}^n$, $K \backslash O$ denotes the subset of elements of $K$ that are not in $O$. By $\mathbb{B}$, we denote the closed unit ball in $\mathbb{R}^n$ centered at the origin. For a function $f : \mathbb{R}^m \rightarrow \mathbb{R}^n$ and a set $D \subset \mathbb{R}^m$, $f(D)$ = $\{ f(x) : x \in D \}$. For a continuously differentiable function $B: \mathbb{R}^n \rightarrow \mathbb{R}$, $\nabla B(x)$ denotes the gradient of the function $B$ evaluated at $x$. By $\mathcal{C}^1$, we denote the set of continuously differentiable functions. Finally, $ F : \mathbb{R}^n \rightrightarrows \mathbb{R}^n $ denotes a set-valued map associating each element $x \in \mathbb{R}^n$ into a subset $F(x) \subset \mathbb{R}^n$.

\section{Preliminaries and Basic Conditions} \label{Sec.1}

\subsection{Hybrid Inclusions}
We consider hybrid systems modeled by 
\begin{align} \label{eq.hsys}
\mathcal{H}: & \left\{ \begin{matrix}  x \in C & ~~~~\dot{x} \in F(x)  \\  x \in D &~~ x^+ \in G(x), \end{matrix} \right.
\end{align}
with the state variable $x \in \mathbb{R}^n$, the flow set $C \subset \mathbb{R}^n$, the jump set $D \subset \mathbb{R}^n$, the flow and the jump set-valued maps, respectively, $F: \mathbb{R}^n \rightrightarrows \mathbb{R}^n$, $G: \mathbb{R}^n \rightrightarrows \mathbb{R}^n $. A solution $x$ to $\mathcal{H}$ is defined on a hybrid time domain denoted $ \dom x \subset \mathbb{R}_{\geq 0} \times \mathbb{N}$. The solution $x$ is parametrized by the ordinary time variable $t \in \mathbb{R}_{\geq 0}$ and the discrete jump variable $j \in \mathbb{N}$. 
Its domain of definition $\dom x$ is such that for each $(T,J) \in \dom x$, $\dom x \cap \left( [0,T] \times \left\{ 0,1,\ldots,J\right\}  \right) = \cup^J_{j=0}\left([t_j,t_{j+1}],j \right)$ for a sequence $\left\{ t_j \right\}^{J+1}_{j=0}$, such that $t_{j+1} \geq t_j$ for each $j \in \left\{0,1, \ldots, J \right\}$ and $t_0 = 0$; see \cite{goebel2009hybrid}. 

A solution $x$ to $\mathcal{H}$, as defined in \ifitsdraft Definition \ref{defsol} 
\else \cite[Definition 12]{draftautomatica}\fi, starting from $x_o$ is said to be complete if it is defined on an unbounded hybrid time domain; that is, the set $\dom x$ is unbounded. Furthermore, it is said to be maximal if there is no solution $y$ to $\mathcal{H}$ such that $ x(t,j) =y(t,j) $ for all $ (t,j) \in \dom x $ with $\dom x$ a proper subset of $\dom y$. Finally, it is said to be nontrivial if $\dom x$ includes at least two points. To simplify notation, we use $x$ to denote the state and also a solution.  We write ``solution $x$'' to make this distinction when needed.

\subsection{Anatomy of Sets}

Different types of cones have been used in the study of differential inclusions. In the following, for a set $K \subset \mathbb{R}^n$, we recall the ones we use in this paper, see \cite{Aubin:1991:VT:120830} for more details. 
\ifitsdraft
\begin{definition}
The \textit{contingent} cone of $K$ at $x$ is given by
\begin{align} \label{eq.toncon} 
T_K(x) := \left\{ v \in \mathbb{R}^n: \liminf_{h \rightarrow 0^+} \frac{|x + h v|_K}{h} = 0 \right\}.
\end{align}
\end{definition}
We also recall the equivalence \cite[Page 122]{aubin2009set}
\begin{align} \label{eq.conti}
 v \in & T_{K}(x) \Leftrightarrow  \nonumber \\ \exists & \left\{ h_i \right\}_{i \in \mathbb{N}} \rightarrow 0^+ ~\mbox{and}~ \left\{v_i\right\}_{i \in \mathbb{N}} \rightarrow v : x + h_i v_i \in K~~\forall i \in \mathbb{N}.
\end{align}
\begin{definition}
The \textit{Dubovitsky-Miliutin} cone of $K$ at $x$ is given by
\begin{align} \label{eq.cone2}
D_K(x):= \left\{v \in \mathbb{R}^n:\exists \epsilon,\alpha >0:x+(0,\alpha](v+\epsilon \mathbb{B}) \subset K \right\}.
\end{align}
\end{definition}
We, also, recall from \cite{Aubin:1991:VT:120830} the following useful property
\begin{align} \label{eq.cone3}
D_K(x)= & \mathbb{R}^n \backslash T_{\mathbb{R}^n \backslash K}(x) = T_{K}(x) \backslash T_{\mathbb{R}^n \backslash K}(x)~~~\forall x \in \partial K.
\end{align}
\begin{definition}
The \textit{external contingent} cone of $K$ at $x$ is given by
\begin{align} \label{eq.exter}
E_K(x) := & \left\{ v\in \mathbb{R}^n :\liminf_{h \rightarrow 0^+} \frac{|x+hv|_K - |x|_K}{h} \leq 0  \right\}. 
\end{align}
\end{definition} 
\else
The \textit{contingent} cone of $K$ at $x$ is given by
\begin{align} \label{eq.toncon} 
T_K(x) := \left\{ v \in \mathbb{R}^n: \liminf_{h \rightarrow 0^+} \frac{|x + h v|_K}{h} = 0 \right\}.
\end{align}
We also recall the equivalence \cite[Page 122]{aubin2009set}
\begin{align} \label{eq.conti}
 v \in T_{K}(x) \Leftrightarrow ~ & \exists \left\{ h_i \right\}_{i \in \mathbb{N}} \rightarrow 0^+ ~\mbox{and}~ \left\{v_i\right\}_{i \in \mathbb{N}} \rightarrow v : \nonumber \\ & x + h_i v_i \in K~~\forall i \in \mathbb{N}.
\end{align}
The \textit{external contingent} cone of $K$ at $x$ is given by
\begin{align} \label{eq.exter}
\hspace{-0.4cm} E_K(x) := \left\{ v\in \mathbb{R}^n :\liminf_{h \rightarrow 0^+} \frac{|x+hv|_K - |x|_K}{h} \leq 0  \right\}. 
\end{align}
The \textit{Dubovitsky-Miliutin} cone of $K$ at $x$ is given by
\begin{align} \label{eq.cone2}
D_K(x) := & \{v \in \mathbb{R}^n : \exists \epsilon >0: \nonumber \\ &
x + \delta (v + w) \in K~\forall \delta \in (0,\epsilon],~\forall w \in \epsilon \mathbb{B}  \}.
\end{align}
\fi

\subsection{Basic Assumptions}

The proposed results on forward invariance and contractivity of a set $K \subset C \cup D$ for $\mathcal{H}$ are obtained under the following standing assumption. 

\noindent
\textbf{Standing assumption.} 
The data of the hybrid inclusion $\mathcal{H} = (C,F,G,D)$ is such that the flow map $F$ is outer semicontinuous and locally bounded with nonempty and convex images,  and $G(x)$ is nonempty for all $x \in D$. Furthermore, the set $K$ is closed. 
\hfill $\bullet$ 

We notice that, in addition to these standing assumptions, the hybrid basic conditions in \cite[Chapter 6]{goebel2012hybrid}, which are not imposed here, also require the sets $C$ and $D$ to be closed and the jump map $G$ to be locally bounded.

Before going further, consider the hybrid inclusion $\mathcal{H} =(C,F,D,G)$ and a closed set $K \subset C \cup D$. Starting from $x_o \in K$, if a solution $x$ leaves the set $K$, then it has to be under one of the two following scenarios:

\begin{enumerate}[label={(Sc\arabic*)},leftmargin=*]
\item \label{item:S1} The solution $x$ leaves the set $K$ after a jump. It implies the existence of $(t,j) \in \dom x$ such that $x(t,j) \in K \cap D$ and 
$(t,j+1) \in \dom x$ with $x(t,j+1) \notin K$ and $x(t,j+1) \in G(x(t,j))$. 

\item \label{item:S2} The solution $x$ leaves the set $K$ by flowing. It implies the existence of $t'_2 > t'_1 \geq 0$ and $j' \in \mathbb{N}$ such that $([t'_1,t'_2] \times \left\{j'\right\}) \subset \dom x $ and $ x( (t'_1,t'_2), j') \subset (U(\partial K) \backslash K) \cap C$, with $x(t'_1, j') \in \partial K $ and $x(t'_2, j') \notin K $. 
\end{enumerate}

In fact, when the set $K$ is closed, under \ref{item:S2}, $x(t'_1,j') \in \partial K \cap K$ and since the solution leaves the set $K$, under \ifitsdraft Definition~\ref{solution definition}\else \cite[Definition 12]{draftautomatica}\fi, $x((t'_1,t'_2],j')$ is a subset of $C\backslash K$ for some $t_2' > t_1'$ sufficiently close to $t_1'$. 

\ifitsdraft
\blue{When the set $K \subset C \cup D$ is not closed, the case in \ref{item:S2} is replaced by the following more general scenario:
\begin{enumerate}[label={(Sc3)},leftmargin=*]
\item \label{item:S3} The solution $x$ leaves the set $K$ by flowing. 
It implies the existence of $t'_2 > t'_1 \geq 0$ and 
$j' \in \mathbb{N}$ such that $([t'_1,t'_2] \times \left\{j'\right\}) \subset \dom x$ and one of the following is true:
\begin{enumerate}
\item $ x( (t'_1,t'_2], j') \subset (U(\partial K) \backslash K) \cap C$, with $x(t'_1, j') \in \partial K \cap K$. 
\item $ x( [t'_1,t'_2), j') \subset K$, with $x(t'_2, j') \in 
\mbox{cl}(K) \backslash K$. 
\end{enumerate}
\end{enumerate}}
\fi

\section{Sufficient Conditions for Forward Pre-Invariance and Invariance Using Barrier Functions}  
\label{Sec.2}

Given a hybrid system $\mathcal{H} = (C,F,D,G)$, for a set $K \subset C \cup D$, following \cite{chai2015notions} and \cite[Definition 6.25]{goebel2012hybrid}, we introduce the two following forward invariance notions.    

\begin{definition}[Forward pre-invariance] \label{defpre}
The set $K$ is said to be forward pre-invariant for 
$\mathcal{H}$ if, for each $x_o \in K$ and each maximal solution $x$ starting from $x_o$, $x(t,j) \in K$ for all $(t,j) \in \dom x$. 
\end{definition} 

\begin{definition}
[Forward invariance] \label{def.com}
The set $K$ is said to be forward invariant for $\mathcal{H}$ if it is forward pre-invariant and for each $x_o \in K$, each maximal solution $x$ starting from $x_o$ is complete. 
\end{definition} 

Furthermore, we assume that the set $K$ is defined as points in $C \cup D$ at which multiple scalar functions are simultaneously nonpositive. These scalar functions form a barrier function candidate defining the set $K$. 

\begin{definition} \label{def.bar}
A function 
$B : \mathbb{R}^n \rightarrow \mathbb{R}^m$ is said to be a barrier function candidate defining the set $K$ if\footnote{$B(x) \leq 0$ means that $B_i(x) \leq 0$ for all $i \in 
\left\{ 1,2,\ldots,m \right\}$.}
\begin{align} \label{eqK}
K = \left\{ x \in C \cup D: B(x) \leq 0 \right\},
\end{align}
where $B(x) := [B_1(x)~B_2(x)~\ldots~B_m(x)]^\top$.
\end{definition}  

If $B$ is continuous, the set $K$ is closed relative to $C \cup D$. If, in addition, $C \cup D$ is closed, then $K$ is automatically closed. 

We introduce the following sets that we use in some statements and proofs. For a set $K$ given as in \eqref{eqK}, we define  
\begin{align} 
K_e := & \left\{ x \in \mathbb{R}^n : B(x) \leq 0 \right\}, \label{eqK1} 
\end{align}
and, for each $i \in \left\{ 1,2,\ldots,m \right\}$,
\begin{align}
K_{ei} := & \left\{ x \in \mathbb{R}^n : B_i(x) \leq 0 \right\}, \label{eqK2} \\
M_i := & \left\{ x \in \partial K : B_i(x) = 0 \right\}. \label{eqK3}
\end{align}

It is useful to notice that $K_e = \cap^m_{i=1} K_{ei}$, $K = K_e \cap (C \cup D)$, and that 
$\partial K = \cup^m_{i=1} M_i \cup \left( \partial K \cap \partial (C \cup D) \right)$. Note that in general $M_i \neq \partial K_{ei}$.

\begin{remark}
In existing literature, motivated by barrier methods for optimization \cite{boyd2004convex}, barrier function candidates\footnote{Barrier functions are also called \textit{potential} functions in \cite{tanner2003stable}.} are introduced as scalar functions that are positive and locally bounded on $\mbox{int}(K)$, and approach infinity as their argument converges to $\partial K$ \cite{WILLS20041415}. The difference between the barrier function therein and the one in 
Definition~\ref{def.bar} is that, the one in \cite{WILLS20041415} guarantees that the solutions starting from $\mbox{int}(K)$ cannot reach the boundary $\partial K$, which in turn 
renders $\mbox{int}(K)$ invariant (in the appropriate sense), see \cite{ames2014control} for a more detailed comparison.
\end{remark}

\begin{remark}
In the case of hybrid systems modeled as hybrid automata, the concept of barrier functions is used to conclude forward invariance, or safety in general, in \cite{prajna2007framework, DAI201762, 10.1007/978-3-642-39799-8_17}. According to these references, for a hybrid automata with $m$ operating modes, a closed set $K_q \subset \mathbb{R}^n$, $q \in \left\{1,2, \dots, Q \right\}$, is associated to each mode (typically determined by the logic variable $q$) and must be forward invariant for the state variable (typically denoted $\zeta$) only during the corresponding $q$ mode. The sets $K_q$ can be different for each mode. Furthermore, since each mode is governed by a differential equation and the state variable $\zeta$ is allowed to jump only when the mode switches, a barrier function candidate $B_q$ is associated with each mode and defines the corresponding set $K_q$ as $K$ in Definition \ref{def.bar}. Our approach covers such construction. Indeed, if we model a hybrid automata as a hybrid inclusion $\mathcal{H}$ while incorporating the mode as a new discrete state variable $q \in \left\{ 1,2, \ldots, Q \right\}$, then, in the augmented space $\mathbb{R}^n \times \left\{1,2, \ldots, Q \right\}$, we can define the set $K := \cup^Q_{q=1} \left( K_q \times \left\{ q \right\} \right)$ and the candidate $B(\zeta,q) := B_q(\zeta)$. It is easy to see that the latter scalar candidate $B$ defines the set $K$ according to Definition \ref{def.bar}.
\end{remark}

\begin{remark}
In our approach, we do not restrict the barrier function candidate to be a scalar function. In general, one is interested in considering a forward invariant set $K$ that is given by multiple inequality constraints being satisfied simultaneously. Also, we notice that it is always possible from \eqref{eqK} to construct a scalar barrier function candidate that defines the closed set $K$ according to Definition \ref{def.bar} as 
\begin{align} \label{scalbar} 
\bar{B}(x) := \max_{ i \in \left\{1, 2, \ldots, m \right\} } B_i(x). 
\end{align}
However, by doing so, if the vector function $B$ is $\mathcal{C}^1$, the resulting barrier function candidate $\bar{B}$ is not guaranteed to be $\mathcal{C}^1$ and it can be only continuous. Indeed, at points $x$ where multiple $B_i$'s are equal, if their gradients are not identical, then $\bar{B}$ is not differentiable at those elements. 
\end{remark}

\subsection{Pre-Invariance Under Standing Assumptions} \label{subsec1}

The results we present in this section are extensions to what has been proposed in \cite{prajna2007framework, DAI201762, 10.1007/978-3-642-39799-8_17} for general  hybrid inclusions while handling the possible noncomplete solutions and using multiple barrier functions instead of only a scalar one. For general differential inclusions (the continuous part of a hybrid inclusion), as pointed out in \cite{Aubin:1991:VT:120830}, forward invariance of a set is a property that depends on the system's dynamics outside the set. Therefore, in the following results, our flow conditions concern only a neighborhood of the boundary $\partial K$ relative to the complement of $K$.

\begin{theorem} \label{prop1}
Given a hybrid system $\mathcal{H} = (C,F,D,G)$ and a $\mathcal{C}^1$ barrier function candidate $B$ defining the set $K$ in \eqref{eqK}. Then, the set $K$ is forward pre-invariant for $\mathcal{H}$ if, for each 
$i \in \left\{1,2,\dots,m \right\}$, there exists a neighborhood $U(M_i)$ such that
\begin{align}
\langle \nabla B_i(x), \eta \rangle \leq & ~0~~~\forall x \in (U(M_i) \backslash K_{ei}) \cap C~\mbox{and} \nonumber \\ & ~~~~~~~~~~~~~~~~~ \forall \eta \in F(x) \cap T_C(x),   
\label{eq.1}   \\ 
B(\eta) \leq & ~0~~~\forall \eta \in G(x),~~\forall x \in D \cap K,  \label{eq.2}   \\
G(x) \subset & ~ C \cup D~~~\forall x \in D \cap K,  \label{eq.2rj}
\end{align}
where $K_{ei}$ and $M_i$ are defined in 
\eqref{eqK2}-\eqref{eqK3} and $T_C$ is the contingent cone of the set $C$.
\end{theorem} 

\begin{proof}  
To prove the statement, we proceed by contradiction. Let us assume that \eqref{eq.1} and \eqref{eq.2rj} hold and the set $K$ is not forward pre-invariant for 
$\mathcal{H}$. That is, there exists a maximal solution $x$ starting from $x_o \in K$ that leaves the set $K$ following one of the scenarios \ref{item:S1} and \ref{item:S2}. First, suppose that the solution $x$ leaves the set $K$ after a jump from $K$ to $ \mathbb{R}^n \backslash K$ following the scenario \ref{item:S1}.  This implies, using \eqref{eq.2rj} and the definition of $B$, the existence of $k \in 
\left\{1,2, \dots, m \right\}$ and $(t,j) \in \dom x$ such that $(t,j+1) \in \dom x$ and $B_k(x(t,j+1))>0$ with $x(t,j+1) \in G(x(t,j))$. However, $x(t,j) \in K \cap D$, hence using \eqref{eq.2}, it follows that $B(x(t,j+1)) \leq 0$; in fact $B(\zeta) \leq 0$ for all $\zeta \in G(x(t,j))$. The latter fact yields a contradiction. Next, suppose that the solution $x$ leaves the set $K$ by flowing under scenario \ref{item:S2}. We conclude in this case that there exists $k \in \left\{1,2, \ldots, m \right\}$ such that $x((t'_1,t'_2], j') \subset (U(\partial K_k) \backslash K_k) \cap C $, where $t'_1$, $t'_2$, and $j'$ are as in \ref{item:S2}. Next, since the function $B$ is assumed to be continuously differentiable and the solution $x(\cdot,j')$ is absolutely continuous on the interval $[t'_1,t'_2]$, it follows that $B(x(\cdot,j'))$ is also absolutely continuous on that interval. By integration, it follows that
\begin{align} \label{eq.integ2}
B_k (x(t'_2,j')) & - B_k (x(t'_1,j')) = \nonumber \\ & \int^{t'_2}_{t'_1} 
\langle \nabla B_k(x(t,j')), \dot{x}(t,j') \rangle 
dt >0
\end{align} 
since $B_k(x(t,j'))>0$ for all $t \in [t'_1, t'_2]$ and $B_k(x(t'_1,j')) = 0$. However, 
$x((t'_1,t'_2],j') \subset (U(\partial K) \backslash K) \cap C $ and, using \ifitsdraft Lemma \ref{lem1}\else\cite[Lemma 2]{draftautomatica}\fi, we conclude that $ \dot{x}(t,j') \in T_C(x(t,j')) $ for almost all $t \in [t'_1, t'_2]$. Moreover, using \eqref{eq.1}, we conclude that, for almost all $t \in (t'_1,t'_2)$, $\langle \nabla B_k (x(t,j)), \eta \rangle \leq 0$ for all $\eta \in F(x(t,j)) \cap T_C(x(t,j))$. Hence, $B_k(x(t'_2,j)) - B_k(x(t'_1,j)) \leq 0$. Hence, the contradiction with \eqref{eq.integ2} follows.
\ifitsdraft 
\else
\hfill $\square$
\fi
\end{proof}

\begin{example} \label{exp1}
Consider the hybrid system $\mathcal{H}$ with the data
\begin{align*}
C := & ~ \left\{ x \in \mathbb{R}^2 : x_2 \geq 0,~x_1 \in [-1,1]  \right\}, \\
F(x) := & ~ \begin{bmatrix}  - x_2^2 \\ x_2 x_1 - x_2 ([2,4] -|x|^2) \end{bmatrix}~~~\forall x \in C, \\
D := & ~ \left\{ x \in \mathbb{R}^2 : x_2 \leq 0,~|x| < 1 \right\}, \\
 G(x) := & ~ [0, x_2] \times [0, |x_1|] ~~~ \forall x \in D.
\end{align*} 
We establish forward pre-invariance for the closed set 
$K := \left\{ x \in C \cup D : |x|^2 \leq 1,~x_2 \geq 0 \right\}$ 
using Theorem \ref{prop1}. To this end, we note that the set $K$ can be written as in \eqref{eqK} using the $\mathcal{C}^1$ barrier function candidate $B(x)=[B_1(x)~~B_2(x)]^\top:=[(|x|^2 - 1)~~-x_2]^\top $. Furthermore, $D \cap K = (-1,1) \times \left\{0 \right\}$ and, for each $x \in D \cap K$, $G(x) = [0 ~~ [0,1]|x_1|]^\top \subset C \cup D$; hence, \eqref{eq.2rj} holds. Moreover, for each $x \in K \cap D$ and for each  $\eta \in G(x)$, there exists $\epsilon \in [0,1]$ such that $B(\eta) = [(\epsilon |x|^2-1)~~-\epsilon |x_1|]^\top \subset \mathbb{R}_{\leq 0} \times \mathbb{R}_{\leq 0}$; thus, \eqref{eq.2} holds. Next, we note that the set $\left( U(M_2) \backslash K_{e2} \right) \cap C$ is empty and one can choose $(U(M_1) \backslash K_{e1}) \cap C = \left\{ x \in C : |x| \in (1,2)   \right\}$. Consequently, for each $\eta \in F(x)$, there exists 
$\epsilon \in [0,2]$ such that $\langle \nabla B_1(x) , \eta  \rangle  = - x_2^2 (2 + \epsilon - |x|^2) \leq  0$ for all $x \in U(M_1) \backslash K_{e1}) \cap C$. Hence, \eqref{eq.1} holds and forward pre-invariance for $\mathcal{H}$ of the set $K$ follows. Note that \eqref{eq.1} does not hold on the entirety of $C \backslash K$. 
\end{example}

In the following example, we apply Theorem \ref{prop1} on a hybrid system including explicit logic variables.

\begin{example}[Thermostat] \label{exptermostat}
Consider the hybrid model of the thermostat system proposed in \cite[Example 1.9]{goebel2012hybrid} with $x:=[q~~z]^\top \in \mathbb{R}^2$,
\begin{align*}
C := & ~ \left( \left\{ 0 \right\} \times C_0 \right)  \cup \left( \left\{ 1 \right\} \times C_1 \right), \\
C_0 := & ~ \left\{z \in \mathbb{R} : z \geq z_{min} \right\}, ~C_1 :=  ~ \left\{z \in \mathbb{R} : z \leq z_{max} \right\}, \\
F(x) := & ~[  0 ~~~~~~ - z + z_o + z_{\Delta} q  ]^\top~~~\forall x \in C, \\
D := & ~ \left( \left\{ 0 \right\} \times D_0 \right) \cup \left( \left\{1 \right\} \times D_1 \right), \\
D_0 := & ~\left\{z \in \mathbb{R} : z \leq z_{min} \right\}, 
~ D_1 :=  ~ \left\{z \in \mathbb{R} : z \geq z_{max} \right\},  \\
G(x) := & ~ [ 1-q ~~~~ z ]^\top~~~\forall x \in D,
\end{align*}
where $z$ is the temperature of the room, $z_o$ represents the natural temperature of the room when the heater is not used, $z_{\Delta}$ the capacity of the heater to raise the temperature in the room by always being on, and $q$ the state of the heater, which is $1$ (on) or $0$ (off). We want to keep the temperature between $z_{min}$ and $z_{max}$ satisfying $ z_o < z_{min} < z_{max} < z_o + z_{\Delta}$. Using Theorem \ref{prop1}, we will show that the set $ K:= \left\{ [q~~z]^\top \in C \cup D: z \in [z_{min}, z_{max}] \right\} = \left\{0,1\right\} \times [z_{min}, z_{max}]$ is forward pre-invariant. To do so, we propose the barrier function candidate $B(x) = [B_1(x)~~B_2(x)]^\top := [ z-z_{max} ~~~~  z_{min}-z ]^\top$. To verify \eqref{eq.2rj}, we note that $C \cup D = \left\{0,1\right\} \times \mathbb{R}$ and $[1-q~~z]^\top \in C \cup D$ for all $[q~~z]^\top \in C \cup D$. Hence, \eqref{eq.2rj} is satisfied. Moreover, $B(G(q,z)) = B([(1-q)~~z]^\top) = B([q~~z]^\top) \leq 0$ for all $[q~~z]^\top \in K \cap D$, the latter inequality holds by definition of the barrier candidate $B$. Hence, \eqref{eq.2} is also satisfied. Finally, to verify \eqref{eq.1}, we note that 
$K_{e1} = \mathbb{R} \times \left(-\infty, z_{max} \right]$, 
$K_{e2} =  \mathbb{R} \times \left[z_{min}, +\infty \right)$, and 
$M_i = \partial K_{ei} \cap (C \cup D)$ for all $i \in \left\{1,2\right\}$. Furthermore, for some $\epsilon>0$, 
$ (U(M_1) \backslash K_{e1}) \cap C = \left\{ 0 \right\} \times (z_{max}, z_{max} +\epsilon)$, and
$ (U(M_2) \backslash K_{e2}) \cap C = \left\{ 1 \right\} \times (z_{min} - \epsilon, z_{min})$. As a result, 
$\langle \nabla B_1(x), F(x) \rangle = z_o - z \leq 0$ for all $x \in (U(M_1) \backslash K_{e1}) \cap C$,
and $\langle \nabla B_2(x), F(x) \rangle = z-z_o-z_\Delta \leq 0$ for all $x \in (U(M_2) \backslash K_{e2}) \cap C$.
\end{example}

\begin{remark} \label{cor1nw}
When the set $K$ is defined as the zero sub-level set of a scalar barrier function candidate $B$; namely, $K = \left\{ x \in C \cup D : B(x) \leq 0 \right\}$ with $m=1$, condition \eqref{eq.1} in Theorem \ref{prop1} reduces to
\begin{align}
\langle \nabla B(x), \eta \rangle \leq 0 &~~~~\forall x \in ( U(\partial K) \backslash K ) \cap C~\mbox{and} \nonumber \\ & ~~~~~~~~~~~~~~~~~  \forall \eta \in F(x) \cap T_C(x). \label{eq.3bis}
\end{align}
\end{remark}

\begin{example}[Boucing ball] \label{expbb}
Consider the bouncing ball hybrid model $\mathcal{H} = (C,F,D,G)$ with $x \in \mathbb{R}^2$,
\begin{align*}
F(x) := & ~ [x_2~~-\gamma]^\top ~~ \forall x \in C, \\
 C := & ~\left\{x \in \mathbb{R}^2: x_1 > 0,~\mbox{or}~x_1 =0~\mbox{and}~x_2 \geq 0  \right\}, \\ 
G(x) := & ~ [0~~-\lambda x_2]^\top ~~ \forall x \in D, \\
D := & ~ \left\{ x \in \mathbb{R}^2: x_1 = 0,~ x_2 \leq 0 \right\}.
\end{align*} 
The constants $\gamma > 0$ and $\lambda \in [0,1]$ are the gravity acceleration and the restitution coefficient, respectively. Consider the closed set 
$$ K := \left\{ x \in C \cup D : 2 \gamma x_1 + (x_2-1) (x_2+1) \leq 0 \right\}. $$ 
The set $K$ can be seen as the sublevel set where the total energy of the ball is less or equal than $1/2$. Hence, $B(x) := 2 \gamma x_1 + (x_2-1) (x_2+1)$ is a barrier function candidate defining the set $K$ as in Definition \ref{def.bar}. To conclude forward pre-invariance of the set $K$ using Theorem \ref{prop1}, we start noticing that $\langle \nabla B(x), F(x) \rangle = 0$  for all $x \in C$; hence, \eqref{eq.3bis} is satisfied. Moreover, for every $x \in K \cap D$, 
$ B(G(x)) = 2\gamma x_1 + \lambda^2 x_2^2 - 1  \leq 2\gamma x_1 +  x_2^2 - 1 \leq 0 $ since $\lambda \in [0,1]$. Hence, \eqref{eq.2} is satisfied. Finally, \eqref{eq.2rj} is satisfied since $G(D) = \left\{0\right\} \times \mathbb{R}_{\geq 0} \subset C \cup D$.
\end{example}

\ifitsdraft 

\begin{example} \label{exp1nw}
Consider the hybrid system $\mathcal{H} = (C,F,D,G)$ with $C$ and $G$ as in Example \ref{exp1}, and $F$ and $D$ given by
\begin{align*}
F(x) := & \begin{bmatrix}  - x_2 x_1 \\ - (|x|^2 - [0,1]) (|x|^2 - \frac{1}{4}) (2-|x|^2) \end{bmatrix}
~~~\forall x \in C, \\
D := & \left\{ x \in \mathbb{R}^2 : x_2 = 0,~|x| \leq 1 \right\}.
\end{align*} 
We employ Theorem \ref{prop1} to verify forward pre-invariance for the closed set \\
$ K := \left\{ x \in C \cup D : |x|^2 - 1 \leq 0,~x_2 \geq 0 \right\}$. This set admits the scalar $\mathcal{C}^1$ barrier function candidate $B(x) := x_2 (|x|^2 - 1)$. According to Remark \ref{cor1nw}, the set $(U(\partial K) \backslash K) \cap C$ can be chosen as $ (U(\partial K) \backslash K) \cap C = \left\{ x \in C : |x_1|<1,~|x| \in (1, 2) \right\} $. Furthermore, it is easy to show that $\langle \nabla B(x), \eta \rangle \leq 0$ for all $\eta \in F(x)$ provided that $0 \leq |x|^2 \leq 2$; hence \eqref{eq.3bis} holds. Furthermore, for any $x \in D \cap K$, $G(x) \subset C \cup D$; hence, \eqref{eq.2rj} holds. Moreover, for any $x \in K \cap D$ and for any $\eta \in G(x)$, there exists $\alpha \in [0,1]$ such that $\eta = [0~~\alpha |x_1|]^\top$ and $B(\eta) = \alpha |x_1| (\alpha^2 |x_1|^2 - 1)$, however, since $x_1 \leq 1$ for all $x \in K \cap D$, \eqref{eq.2} follows. 
\end{example}

\fi

\begin{remark} \label{rem.compar}
 The flow condition \eqref{eq.1} in Theorem \ref{prop1} is more general than those in \cite{prajna2007framework, ames2014control, glotfelter2017nonsmooth, DAI201762, 10.1007/978-3-642-39799-8_17} in the sense that the inequality in \eqref{eq.1} does not need to hold on the entire set $C$. In \cite{ames2014control, glotfelter2017nonsmooth}, the flow condition \eqref{eq.1} is expressed as
\begin{align*}
\langle \nabla B(x), \eta \rangle \leq \rho(B(x)) & ~~~~ \forall x \in C ~~~~ \forall \eta \in F(x), 
\end{align*}
where $\rho: \mathbb{R} \rightarrow \mathbb{R}$ is an extended class-$\mathcal{K}$ function; 
namely, $\rho(0) = 0$ and $B \mapsto \rho(B)$ is strictly decreasing. Furthermore, in \cite{DAI201762, 10.1007/978-3-642-39799-8_17}, the function 
$\rho: \mathbb{R} \rightarrow \mathbb{R}$ is assumed to be only  locally Lipschitz with $\rho(0) = 0$. Hence, the function $\rho$ is allowed to be positive provided that its growth is bounded locally by a linear function. We consider this type of relaxation in Section \ref{Subsec.3}. This being said, for continuous-time systems with inputs, when the flow inequalities hold on the entire state space, as shown in \cite{jankovic2018robust}, numerical methods can be employed to compute an input value that assures safety.
\end{remark}

\begin{remark} \label{remadded}
From conditions \eqref{eq.1} in Theorem \ref{prop1}, it is straightforward to conclude that it is enough for each barrier function candidate $B_i$ to be of class $\mathcal{C}^1$ only on a neighborhood of the boundary $M_i$. 
\end{remark}

\begin{remark} \label{remadded-}
In Theorem \ref{prop1} (as well as in upcoming results), the jump condition \eqref{eq.2} in Theorem \ref{prop1} can be formulated using a different barrier function candidate than the one used to formulate the flow condition \eqref{eq.1} in Theorem \ref{prop1}. However, the two different barrier function candidates still need to define the same set $K$, according to \eqref{eqK}. In this paper, for simplicity, we present results using the same barrier function candidate in the flow and the jump conditions, but extensions to the case where they are different are straightforward.   
\end{remark}

\subsection{Pre-Invariance Under Lipschitz-Like Flow Map} \label{Subsec.1}

In Theorem \ref{prop1}, the inequality in \eqref{eq.1} needs to be satisfied on a neighborhood outside the set $K$ rather than just on $\partial K$. This is also the case for $m=1$ for which \eqref{eq.3bis} is required; see Remark \ref{cor1nw}. To assess the possibility of relaxing such requirement, we consider the flow condition
\begin{align} 
\langle \nabla B_i(x), \eta \rangle \leq & ~0~~~\forall x \in M_i \cap C, ~~ \forall \eta \in F(x). \label{eq.9} 
\end{align}
We will show that under some regularity assumptions on $F$ and $B$, condition \eqref{eq.9} can be used to conclude forward pre-invariance of the set $K \subset C \cup D$. Furthermore, as we show in Example \ref{expcount}, when we further relax \eqref{eq.9} to
\begin{align*} 
\langle \nabla B_i(x), \eta \rangle \leq & ~0~~\forall x \in M_i \cap C, ~ \forall \eta \in F(x) \cap T_C(x),  
\end{align*}
we fail to conclude forward pre-invariance of the set $K$.

When \eqref{eq.9} and \eqref{eq.2}-\eqref{eq.2rj} hold, the closed set $K \subset C \cup D$ can fail to be forward pre-invariant for the reasons enumerated below, where we assume without loss of generality that $m=1$ (i.e., $B$ is a scalar function).

\begin{enumerate} [label={\arabic*)},leftmargin=*]
\item Assume the existence of $x_o \in \partial K \cap \mbox{int}(C)$ such that $\nabla B(x_o) = 0$. Assume also that $F(x_o) \subset D_{\mathbb{R}^n \backslash K}(x_o)$, where $D_{\mathbb{R}^n \backslash K}$ is the Dubovitsky-Miliutin cone of the set $\mathbb{R}^n \backslash K$. In this case, condition \eqref{eq.9} is satisfied at $x_o$. However, according to \ifitsdraft Theorem \ref{thm3} \else \cite[Theorem 6]{draftautomatica}\fi, there exists a nontrivial solution starting from $x_o$ and flowing outside the set $K$. Hence, the set $K$ is not forward pre-invariant although \eqref{eq.9} is satisfied; cf. \cite[Proposition 1]{ames2014control}. To handle this situation, one needs to assume that the gradient of $B$ is non-degenerate on 
$\partial K_e \cap C$; namely, 
\begin{align}
\nabla B (x) \neq & 0~~~\forall x \in \partial K_e \cap C.   
 \label{eq.9-1-} 
\end{align} 

\item When the solutions starting from $x_o \in \partial K \cap C$ are nonunique, even if $\nabla B(x_o) \neq 0$, when $\langle \nabla B(x_o), \eta \rangle = 0$ for any $\eta \in F(x_o)$, we can always consider the existence of a solution starting from $x_o \in \partial K \cap C$ with a speed that is tangent to $\partial K$ but leaving the set $K$. Such a scenario is illustrated in Example \ref{expillu} below. We also notice that this pathology does not occur when $\dot{x} \in  F(x)$ for all $x \in C$ has unique solutions.
\end{enumerate}

The following example is inspired by \cite[Page 1751]{Blanchini:1999:SPS:2235754.2236030}.
\begin{example} \label{expillu}
 Consider the two-dimensional differential equation
$ \dot{x} = [ 1 ~~ \sqrt{|x_2|} ]^\top =: F(x)$ and the set $ K:= \left\{ x \in \mathbb{R}^2:~x_2 \leq 0 \right\}$. Note that this system can be interpreted as a hybrid system with $C= \mathbb{R}^2$, $D$ empty, and $G$ arbitrary. The set $K$ can be defined using the barrier function candidate $B(x) := x_2 $ satisfying $|\nabla B (x)| = 1 \neq 0 $ and $\langle \nabla B(x), F(x) \rangle = 0$ for all $x \in \partial K$. However, the set $K$ is not forward pre-invariant since $x(t)=[t~~~(1/4)t^2]^\top$ defined for all $t \geq 0$ is a solution starting from $x(0) = 0 \in K$ that leaves the set $K$. 
\end{example}

The latter example confirms the fact that, for differential inclusions, forward pre-invariance of a closed set $K$ is a property of the system outside the set $K$ rather than on its boundary or in its interior. However, when the flow map $F$ satisfies extra regularity conditions outside the set $K$, it is possible to restrict the conditions in \eqref{eq.1} and \eqref{eq.3bis} to hold only on the boundary $\partial K$. This is possible, for example when $F$ is locally Lipschitz as shown in \cite{nagumo1942lage, bony1969principe, brezis1970characterization} for differential equations and in \cite{Aubin:1991:VT:120830} for differential inclusions. In the aforementioned references, contingent-cone-based conditions are used and shown to be necessary as well as sufficient, provided that the system's dynamics is defined on an open set containing the closed set $K$. For differential equations  defined in $\mathbb{R}^n$ with locally Lipschitz right-hand side, the latter contingent-cone-based conditions are expressed in terms of a scalar barrier function candidate 
in \cite{ames2014control}.

On the other hand, to conclude forward pre-invariance using flow conditions satisfied only on the boundary of the set $K$, the Lipschitz regularity of the flow map $F$ can be relaxed by modifying the right-hand side in 
\ifitsdraft \eqref{eq.lipset} \else \cite[Definition 15]{draftautomatica} \fi using uniqueness functions. The latter is shown in \cite{redheffer1972theorems} for differential equations.

\begin{definition} [Uniqueness function]\label{defuniq}
A function $\rho : \mathbb{R} \rightarrow \mathbb{R}$ is said to be a uniqueness function if, for each continuous function 
$\delta: \mathbb{R}_{\geq 0} \rightarrow \mathbb{R}_{\geq 0}$ such that $\delta(0) = 0$ and there exists $\epsilon > 0$ such that
\begin{align} \label{equniqad} 
\hspace{-0.2cm} \limsup_{h \rightarrow 0^+} \frac{\delta(t+h)-\delta(t)}{h} \leq & ~\rho(\delta(t))~
\mbox{for a.a.}~t \in [0, \epsilon],
\end{align} 
it follows that $\delta(t) = 0$ for all 
$t \in [0,\epsilon]$.
\end{definition}

\begin{remark}
Uniqueness functions, in this paper, are introduced in a slightly different way compared to \cite{redheffer1972theorems}. Indeed, in this reference, a function $\rho : \mathbb{R} \rightarrow \mathbb{R}$ is said to be a uniqueness function if, for each continuous function $\delta: \mathbb{R}_{\geq 0} \rightarrow \mathbb{R}_{\geq 0}$ such that $\delta(0) = 0$ and there exists $\epsilon > 0$ such that, in addition to \eqref{equniqad}, we have
\begin{align*}
\limsup_{h \rightarrow 0^+} \frac{\delta(t)-\delta(t-h)}{h} \leq & ~\rho(\delta(t))~~
\mbox{for a.a.}~t \in [0, \epsilon],
\end{align*}
it follows that $\delta (t) = 0$ for all $t \in [0,\epsilon]$. That is, the notion of uniqueness function used in this paper considers only the upper-right Dini derivative.
\end{remark}

In the following statement, we propose flow and jump conditions that are sufficient and need to be satisfied only on elements of the set $K$ provided that the following assumption holds. 

\begin{assumption} \label{ass1}
For every $x \in \partial K_e \cap C$,
\begin{equation} \label{eq.9-1}
\begin{aligned} 
\exists v \in \mathbb{R}^n : \langle \nabla B_i(x), v \rangle < 0 \\ \forall i \in \left\{1,2,\ldots,m\right\} &~ \mbox{s.t.} ~ B_i(x)=0.     
\end{aligned} 
\end{equation}
\end{assumption}

Assumption \ref{ass1} is known as \textit{transversality condition} in \cite{aubin2009set} and allows to define the contingent cone $T_K$ at the intersection between different zero-level sets --- see \ifitsdraft Lemma \ref{lem3} below \else \cite[Lemma 3]{draftautomatica}\fi. Furthermore, Assumption \ref{ass1} reduces to \eqref{eq.9-1-} when the barrier function $B$ is scalar. 

\begin{theorem} \label{propnou}
Given a hybrid system $\mathcal{H} = (C,F,D,G)$ and a 
$\mathcal{C}^1$ barrier function candidate $B$ defining the set $K$ in \eqref{eqK}. Furthermore, assume that there exists a neighborhood $U(\partial K)$ such that, for each $x \in ( U(\partial K) \backslash K ) \cap C $ and for each $y \in \partial (K \cap C)$,
\begin{align} 
\hspace{-0.3cm} (x-y)^\top F(x) \subset (x-y)^\top F(y) + |x-y| \rho(|x-y|) \mathbb{B}  \label{eq.07prop}
\end{align}
with $\rho : \mathbb{R} \rightarrow \mathbb{R}$ a uniqueness function. Then, the closed set $K$ is forward pre-invariant provided that \eqref{eq.2}-\eqref{eq.2rj}, \eqref{eq.9}, and Assumption \ref{ass1} hold, and there exists a neighborhood 
$U(\partial K_e \cap \partial C)$ such that one of the following extra conditions holds:
\begin{enumerate}[label={(\alph*)},leftmargin=*]
\item \label{item:a} For any $x \in \left( U(\partial K_e \cap \partial C) \cap \partial K_e \right) \backslash C$, \eqref{eq.9-1} holds and
\begin{equation}
\label{eqraj1}
\begin{aligned} 
\hspace{-0.4cm} \langle \nabla B_i(x), \eta \rangle & \leq 0~~\forall \eta \in F(x)~\mbox{and} \\ & \forall i \in \left\{ 1, 2,\ldots, m \right\}~ \mbox{s.t.}~ B_i(x) = 0.
\end{aligned}
\end{equation}

\item \label{item:b} 
For any $x \in U(\partial K_e \cap \partial C) \cap \partial K \cap \partial C$,
\begin{align}
F(x) \subset T_{K \cap C}(x). \label{eqraj2}
\end{align} 
\item \label{item:c} The set $C$ is convex and \eqref{eqraj2} holds for all $ x \in \partial K_e \cap \partial C $.
\end{enumerate}
\end{theorem}

\begin{proof}
The proof that the solutions starting from the set $K$ cannot jump outside according to \ref{item:S1}, under \eqref{eq.2}-\eqref{eq.2rj}, is the same as in the proof of Theorem \ref{prop1}. Next, we prove that the solutions starting from $\partial K_e \cap C$ cannot leave the set $K$ by flowing as in scenario \ref{item:S2}. For this purpose, we adapt the steps of the proof presented in \cite[Proof of Theorem 1]{redheffer1972theorems} to our more general setting. Let $t'_1 \geq 0$ and $t'_2 > t'_1$ be such that there exits a solution flowing from $x_o := x(t'_1, 0) \in \partial K \cap C$ and satisfying $x(t,0) \in ( U(\partial K) \backslash K ) \cap C$ for all $t \in (t'_1,t'_2)$. We use $y_1$ to denote the projection of $x(t,0)$ on the set $K_e$ and $y_2$ to denote the projection of $x(t,0)$ on the set $K \cap C$. Furthermore, we define $\delta_1(t) := |x(t,0)-y_1|$ and $\delta_2(t) := |x(t,0)-y_2|$. It follows that $\delta_i(t'_1) = 0$ and $\delta_i(t) > 0$ for all $t\in (t'_1, t'_2)$, since $x(t,0) \in (U(\partial K) \backslash K) \cap C$ for all $t\in (t'_1, t'_2)$,  for each $i \in \left\{1,2\right\}$. Using the identity $ a-b = (a^2 - b^2)/(a+b) $ for $a$ and $b$ nonnegative and for any $h>0$ such that $t$ and $t+h$ in $(t'_1,t'_2)$, we derive the inequality
\begin{align} \label{eq.ns1}
\delta_i(t+h) - & \delta_i(t) \leq |x(t+h,0)-y_i| - |x(t,0) - y_i| \nonumber \\ & = \frac{|x(t+h,0)-y_i|^2 - |x(t,0) - y_i|^2}{|x(t+h,0) - y_i| + |x(t,0) - y_i|}
\end{align}
for each $i \in \left\{1,2\right\}$. Furthermore, for almost all $t \in (t'_1,t'_2)$, we replace $x(t+h,0)$ by 
\begin{align} 
x(t+h,0) = x(t,0) + h \dot{x}(t,0) + o(h)
\end{align} 
with $o(h)$ the remainder of the first order Taylor expansion of $h \mapsto x(t+h)$ around $h=0$, which satisfies
$\lim_{h \rightarrow 0} o(h)/h = 0$. Using the previous limit and the inequality 
$$ |x(t+h,0) - y_i| \leq |x(t,0) - y_i| + h |\dot{x}(t,0) + o(h)/h|, $$ 
we obtain that, for all $i \in \left\{1,2\right\}$ and for almost all $t \in (t'_1, t'_2)$,
\begin{align} \label{eq.ns2} 
\limsup_{h \rightarrow 0^+} \frac{\delta_i(t+h) - \delta_i(t)}{h} \leq \frac{(x(t,0)-y_i)^\top \dot{x}(t,0)}{|x(t,0)-y_i|}. 
\end{align}

Next, we have the following claim: 
\begin{claim} \label{clm1}
Under \eqref{eq.9} and \eqref{eq.9-1}, the following is true.
\begin{enumerate}[label={(cl\arabic*)},leftmargin=*]
\item \label{item:cl1} If \ref{item:a} holds, then $(x(t,0)-y_1)^\top \eta_y \leq 0$ for all $\eta_y \in F(y_1)$.
\item \label{item:cl2} If \ref{item:b} or \ref{item:c} hold, then $(x(t,0)-y_2)^\top \eta_y \leq 0$ for all $\eta_y \in F(y_2)$.
\end{enumerate}
\end{claim}
To prove the claim, we proceed as follows:
\begin{itemize}
\item Under \ref{item:a} in Theorem \ref{propnou}, 
we show that $F(y_1) \subset T_{K_e}(y_1)$. Indeed, under \eqref{eq.9} and \eqref{eq.9-1} and using \ifitsdraft  Lemma \ref{lem3}\else \cite[Lemma 3]{draftautomatica}\fi, it follows that $F(y_1) \subset T_{K_e}(y_1)$ when $y_1 \in \partial K_e \cap C$. Similarly, using the same argument under \ref{item:a}, when $y_1 \in \left( U(\partial K_e \cap \partial C) \cap \partial K_e \right) \backslash C$, we also have $F(y_1) \subset T_{K_e}(y_1)$. Finally,  we use the fact that when $ t'_2$ is sufficiently small, the element $y_1$ corresponding to the projection of $x(t,0) \in ( U(\partial K) \backslash K ) \cap C$ on the set $K_e$ belongs necessarily to either $ \partial K_e \cap C $ or $\left( U(\partial K_e \cap \partial C) \cap \partial K_e \right) \backslash C$. 

\item Under \ref{item:b} or \ref{item:c} in Theorem \ref{propnou}, we show that $F(y_2) \subset T_{K \cap C}(y_2)$. Indeed, under \eqref{eq.9}, Assumption \ref{ass1}, and using \ifitsdraft  Lemma \ref{lem3}\else \cite[Lemma 3]{draftautomatica}\fi, it follows that $F(y_2) \subset T_{K_e}(y_2) = T_{K \cap C}(y_2)$ when $y_2 \in \partial K_e \cap \mbox{int}(C)$. Next, under \ref{item:b}, when $y_2 \in U(\partial K_e \cap \partial C) \cap \partial K \cap \partial C$, we also have $F(y_2) \subset T_{K \cap C}(y_2)$. Finally, we use the fact that when $x_o \in \partial K_e \cap C$ and for $ t'_2 >0$ sufficiently small, the element $y_2$ corresponding to the projection of $x(t,0)$ on the set $K \cap C$ belongs necessarily to either $ (\partial K_e \cap \mbox{int}(C))$ or $(U(\partial K_e \cap \partial C) \cap \partial K \cap \partial C$. Furthermore, under \ref{item:c}, the set $C$ is convex and we show that $y_2 \in \partial K_e \cap C$. Hence, it becomes enough to have $F(y_2) \subset T_{K \cap C}(y_2)$ when $y_2 \in \partial K_e \cap \partial C$, which is true under \ref{item:c}. Finally, to show that $y_2 \in \partial K_e \cap C$ when $C$ is convex, we use the fact that $y_2$ corresponds to a projection on the set $K \cap C$. Hence, either $ y_2 \in \partial K_e \cap C$ or $y_2 \in \partial C \cap \mbox{int}(K_e)$. We propose to exclude the latter case using contradiction. That is, assume that $y_2 \in \partial C \cap \mbox{int}(K_e)$ and consider the line segment relating $y_2$ to $x(t,0)$ denoted by $[y_2, x(t,0)]$. Since both $y_2$ and $x(t,0)$ lie in the set $C$ and since the set $C$ is convex it follows that the segment $[y_2,x(t,0)]$ also belongs to $C$. Furthermore, since $x(t,0) \in C \backslash K_e$ and $y_2 \in \mbox{int}(K_e)$ it follows the existence of $y_o$ belonging to the open segment $(x(t,0), y_2)$ such that $y_o \in \partial K_e \cap C = \partial K \cap C$. Hence $|x(t,0) - y_o| < |x(t,0) - y_2| = \min_{ z \in K \cap C } \left\{ |x(t,0) - z| \right\}$, which yields to a contradiction.      
\end{itemize}
Now, to conclude \ref{item:cl1} and \ref{item:cl2}, we use the inequalities 
\begin{align*} 
 |x(t,0)-y_1| & \leq |x(t,0)- y_1 - h \eta_y| + |y_1 + h \eta_{y1}|_{K_e}, \\
 |x(t,0)-y_2| & \leq |x(t,0)-y_2-h \eta_y| + |y_2 + h \eta_{y2}|_{K \cap C},
\end{align*}
where $(\eta_{y1}, \eta_{y2}) \in F(y_1) \times F(y_2)$. To obtain the previous inequalities, we used the fact that the functions $|\cdot|_{K_e}$ and $|\cdot|_{K \cap C}$ are globally Lipschitz with a Lipschitz constant equal to $1$, $|x(t,0)-y_1| = |x(t,0)|_{K_e}$, and $|x(t,0)-y_2| = |x(t,0)|_{K \cap C}$. Next, by taking the square in both sides of the last two inequalities, and dividing by $h$, we obtain, for each $\eta_{y_1} \in F(y_1)$,
\begin{align*}
|x(t,0) - y_1|^2 & /h  \leq |x(t,0)-y_1-h \eta_{y_1}|^2/h + 
\nonumber \\ & |y_1 + h \eta_{y_1}|_{K_e}^2/h + 
\nonumber \\ & 2 |x(t,0) - y_1 - h \eta_{y_1}| |y_1 + h \eta_{y_1}|_{K_e}/h, 
\end{align*}
which implies that
\begin{align*} 
|x(t,0) & - y_1|^2 /h \leq |x(t,0)-y_1|^2/h + h |\eta_{y_1}|^2 - 
\nonumber \\ & 2 (x(t,0)-y_1) \eta_{y_1} + h (|y_1 + h \eta_{y_1}|_{K_e}/h)^2 +
\nonumber \\ & 2 |x(t,0)-y_1-h \eta_{y_1}| 
|y_1+h \eta_{y_1}|_{K_e}/h. 
\end{align*}
Similarly, for each $\eta_{y_2} \in F(y_2)$,
\begin{align*}
|x(t,0)-y_2|^2 & /h \leq |x(t,0)-y_2|^2/h + h |\eta_{y_2}|^2 - 
\nonumber \\ & 2 (x(t,0)-y_2) \eta_{y_2} +
\nonumber \\ & h (|y_2 + h \eta_{y_2}|_{K \cap C}/h)^2 +
\nonumber \\ & 2 |x(t,0)-y_2-h \eta_{y_2}| 
|y_2 + h \eta_{y_2}|_{K \cap C}/h. 
\end{align*}

Finally, letting $h \rightarrow 0^+$ through a suitable sequence, \ref{item:cl1} and \ref{item:cl2} are proved using the fact that $ \liminf_{h \rightarrow 0^+} |y_1 + h \eta_{y_1}|_{K_e} / h = 0$ and $\liminf_{h \rightarrow 0^+} |y_2 + h \eta_{y_2}|_{K \cap C}/h = 0 $ since we already showed that $\eta_{y_1} \in T_{K_e}(y_1)$ under \ref{item:a} in Theorem \ref{propnou} and $\eta_{y_2} \in T_{K \cap C}(y_2)$ under \ref{item:b} or \ref{item:c} in Theorem \ref{propnou}. 

Using the Claim \ref{clm1}, for each $\eta_y \in F(y_i)$, the term 
$ - \frac{(x(t,0) - y_i)^\top \eta_y}{|x(t,0)-y_i|}$ can be added in \eqref{eq.ns2} which then can be rewritten as
\begin{align} \label{eq.ns3} 
\limsup_{h \rightarrow 0^+} \frac{\delta_i(t+h) - \delta_i(t)}{h} \leq & \frac{(x(t,0)-y_i)^\top (\dot{x}(t,0) - \eta_y)}{|x(t,0)-y_i|},
\end{align}
where $i = 1$ if \ref{item:a} in Theorem \ref{propnou} holds and $i = 2$ if \ref{item:b} or \ref{item:c} in Theorem \ref{propnou} hold. Since $\dot{x}(t,0) \in F(x(t,0))$ for almost all $t \in (t'_1, t'_2)$, using \eqref{eq.07prop} with $x = x(t,0)$ and $y = y_i$, it is always possible to find $\eta^*_y \in F(y_i)$ such that
$$ (x(t,0) - y_i)^\top (\dot{x}(t,0) - \eta^*_y)^\top \leq |x(t,0) - y_i| \rho (|x(t,0) - y_i|). $$ 
Applying this inequality to \eqref{eq.ns3} and replacing $\eta_y$ therein by $\eta^*_y$, we obtain 
\begin{align*} 
\limsup_{h \rightarrow 0^+} \frac{\delta_i(t+h) - \delta_i(t)}{h} \leq & \rho(|x(t,0)-y_i|) = \rho(\delta_i(t)) 
\end{align*}
for almost all $t \in (t'_1, t'_2)$ with $i = 1$ if \ref{item:a} holds and $i = 2$ \ref{item:b} or \ref{item:c} hold. Since the function $\rho$ is a uniqueness function, it follows that 
$\delta_i (t) = 0$ for all $t \in (t'_1,t'_2)$. Hence, the contradiction follows. 
\ifitsdraft 
\else
\hfill $\square$
\fi 
\end{proof}

\begin{remark}
Note that condition \eqref{eq.07prop} holds for free
when $F$ is locally Lipschitz. In fact, every locally Lipschitz map $F$ satisfies \eqref{eq.07prop} with $\rho (\omega) = k \omega$ for some $k > 0$. Condition \eqref{eq.07prop} is more general than Lipschitzness and allows for functions $\rho$ that are not necessarily linear. In particular, when $F$ is such that \eqref{eq.07prop} holds for $\rho(\omega) := \omega \log \omega$, then $\rho$ is a uniqueness function \cite{agarwal1993uniqueness}. The latter function belongs to the more general class of \textit{Osgood} functions that are uniqueness functions but not necessarily locally Lipschitz \cite{redheffer1972theorems}.
\end{remark}

\begin{remark}
The flow condition \eqref{eqraj1} in Theorem \ref{propnou} is a reinterpretation, in terms of barrier functions, of the well-known contingent cone-based condition used in \cite{nagumo1942lage, bony1969principe, brezis1970characterization, Aubin:1991:VT:120830} provided that Assumption \ref{ass1} holds. 
\end{remark}

\ifitsdraft
\blue{
\begin{remark}
In \cite{ames2014control}, when $\mathcal{C} = \mathcal{D}$ therein, the proposed result for differential equations is a particular case of Theorem \ref{propnou}. Indeed, consider the system 
$$ \dot{x} = f(x)~~~x \in \mathbb{R}^n, $$
where $f: \mathbb{R}^n \rightarrow \mathbb{R}^n$ is locally Lipschitz. Consider a closed set $K \subset \mathbb{R}^n$ admitting a scalar and continuously differentiable barrier function candidate $B$. According to Theorem \ref{propnou}, when $\nabla B(x) \neq 0$ for all $x \in \partial K$, the set $K$ is forward pre-invariant if 
\begin{align} 
\langle \nabla B (x), f(x) \rangle \leq 0 ~~~~  \forall x \in \partial K. \label{eq.2c2}
\end{align}
In \cite{ames2014control}, \eqref{eq.2c2} is replaced by
\begin{align} 
\langle \nabla B (x), f(x) \rangle \leq \gamma(B(x)) ~~~~  \forall x \in K, \label{eq.2c2bis}
\end{align}
where $\gamma$ is a class-$\mathcal{K}$ function; thus, 
$\gamma(0) = 0$. The latter fact implies that \eqref{eq.2c2} is tighter than \eqref{eq.2c2bis}. 
\end{remark}}
\fi

In Example \ref{expcount}, we show that when none of the conditions \ref{item:a}-\ref{item:c} in Theorem \ref{propnou} is satisfied, there exist situations where $K$ fails to be forward pre-invariant. Also, we show that the aforementioned conditions are only sufficient.

\begin{example} \label{expcount}
Consider the differential inclusion 
$\mathcal{H} = (C,F,\emptyset,\star)$, where 
$F(x) := [1 ~~ 0]^\top$ for all $x \in \mathbb{R}^2$ and
\begin{align*} 
C  := & \left\{ x \in \mathbb{R}^2 : |x_2| \geq x_1^2 \right\}  \cup \left\{ x \in \mathbb{R}^2 : x_1 \leq 0 \right\} \\ & \cup \left\{ x \in \mathbb{R}^2 : x_2 = 0 \right\}.
\end{align*}
 Furthermore, consider the barrier function candidate 
\begin{align*}
B(x) := \left\{ \begin{matrix} x_2 & ~~ \mbox{if}~x_1 \leq 0 \\  x_2+ x_1^3  &  ~~ \mbox{otherwise} \end{matrix}  
\right.
\end{align*}
defining a closed set $K$ according to \eqref{eqK}. 
\ifitsdraft
\else
\begin{figure}
\begin{center}
\def\svgwidth{ \columnwidth}
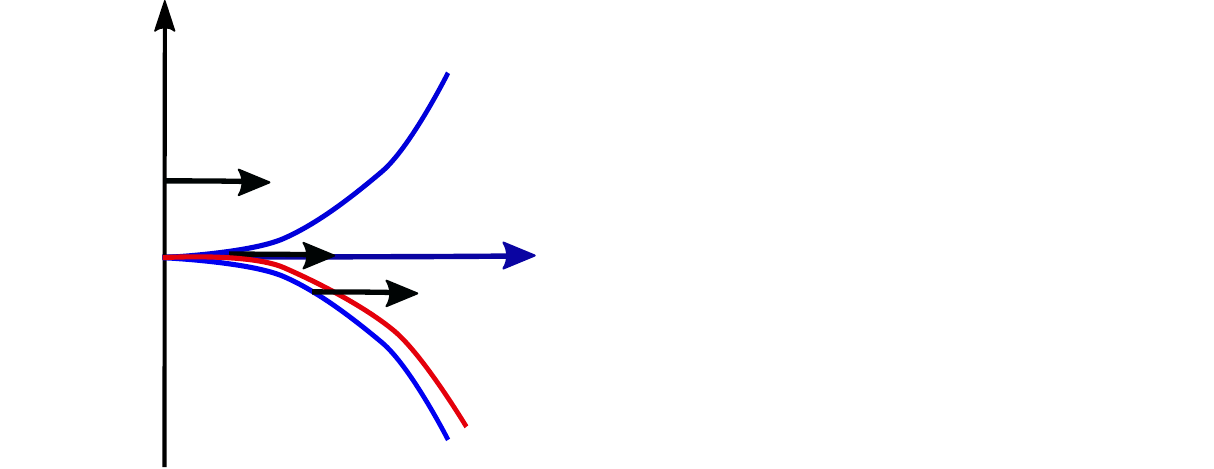
\end{center}
\caption{Illustration of the system in Example \ref{expcount}.}
\end{figure}
\fi
We will show that \eqref{eq.9}-\eqref{eq.9-1-} are satisfied; however, none of the conditions \ref{item:a}-\ref{item:c} in Theorem \ref{propnou} is satisfied. As a result, we show that the set $K$ is not forward pre-invariant. Indeed, note that 
$\nabla B(x) = [0~~1]^\top$ if $x_1 \leq 0$ and 
$ \nabla B(x) = [2 x_1^2~~1]^\top$ otherwise; hence, \eqref{eq.9-1-} is satisfied. Furthermore, for each 
$x \in \partial K_e \cap C = \left\{ x \in \mathbb{R}^2 : x_2 = 0,~~x_1 \leq 0 \right\}$, we have $\langle \nabla B(x), F(x) \rangle = \langle [0~~1]^\top, [1~~0]^\top \rangle \leq 0 $; hence, \eqref{eq.9} is satisfied. Moreover, \eqref{eq.2}-\eqref{eq.2rj} are trivially satisfied since $D = \emptyset$. Next, it is easy to see that the set $C$ is not 
convex; hence, condition \ref{item:c} is not satisfied. Furthermore, when $x \in \left( U(\partial K_e \cap \partial C) \cap \partial K_e \right) \backslash C = \left\{ x \in \mathbb{R}^2 : x_2 = - x_1^3,~ x_1 \in (0, \epsilon],~\epsilon > 0 \right\}$, $\langle \nabla B(x), F(x) \rangle = \langle [2 x_1^2~~1]^\top, [1~~0]^\top \rangle = 2 x_1^2 > 0$; hence, \ref{item:a} is not satisfied. Next, for any $x \in U(\partial K_e \cap \partial C) \cap \partial K \cap \partial C = \left\{ x \in \mathbb{R}^2 : x_2 = - x_1^2,~ x_1 \in [0, \epsilon],~\epsilon > 0 \right\}$, if $x = 0$, then $F(0) \in T_{K \cap C}(0)$; however, if $x \neq 0$, then $F(x) \notin T_C(x)$ because $F(x)$ points outside the set $\left\{ x \in \mathbb{R}^2 : |x_2| \geq x_1^2 \right\}$, which defines $C$ when $x_1>0$. Hence, \ref{item:b} is also not satisfied. Finally, the constrained differential inclusion $\mathcal{H} = (C, F, \emptyset, \star)$ admits the solution $x(t) = [t~~0]^\top$, $t \geq 0$, starting from $x_o = 0 \in K$, that leaves the set $K$.  

Now, in order to show that none of the conditions \ref{item:a}-\ref{item:c} is necessary, we slightly modify the set $C$ in order to render the set $K$ forward pre-invariant while maintaining \ref{item:a}-\ref{item:c} unsatisfied. That is, consider the new flow set 
$$ C_1 := \left\{ x \in \mathbb{R}^2 : |x_2| \geq x_1^2 \right\} \cup \left\{ x \in \mathbb{R}^2 : x_1 \leq 0 \right\}. $$
Note that the unique solution starting from each $x_o := [x_{o1}~~x_{o2}]^\top \in K$ satisfies $x(t) = [x_{o1} + t~~x_{o2}]^\top$ for all $t \geq 0$ and cannot leave the set $K$ while remaining in $C$. However, using the same arguments as in the previous paragraph, we conclude that conditions \ref{item:a}-\ref{item:c} remain unsatisfied.  
\end{example}

\begin{example}  
Consider the hybrid system $\mathcal{H}$ with the data
$C := \left\{ x \in \mathbb{R}^2 : |x_1| \leq 1 \right\},~F(x) := \begin{bmatrix} 1 \\ - [0, x_2 \log |x_2|] \end{bmatrix}$, \\
$D := \left\{ x \in \mathbb{R}^2 : x_1 = 1 \right\},~~~G(x) := [-1 ~~~~ x_2]^\top$. 
\\ 
To conclude forward pre-invariance of the set $K := \left\{ x \in C \cup D : x_2 \leq 0 \right\}$ admitting the 
$\mathcal{C}^1$ barrier function candidate 
$B(x) := x_2$, we start noting that conditions \eqref{eq.2}-\eqref{eq.2rj} are satisfied since $G(x) \in C$ for all $x \in D$ and $B(G(x)) = B(x) \leq 0$ for all $x \in K$. Furthermore, conditions \eqref{eq.9}-\eqref{eq.9-1-} are also satisfied since $\nabla B(x) =[0~~1]^\top$ for all $x \in \mathbb{R}^2$ and 
$ \langle \nabla B(x), F(x) \rangle =  - [0, x_2 \log |x_2|] \leq 0$ for all $x \in \partial K_e$. For the same reason, condition \ref{item:a} in Theorem \ref{propnou} is also satisfied. Finally, it remains to show that for any $x \in ( U(\partial K) \backslash K ) \cap C $ and for any $y \in \partial (K \cap C)$, \eqref{eq.07prop} is satisfied. Indeed, for $y := [y_1~y_2]^\top \in \partial K \cap C$, it follows that $y_2 = 0$. At the same time, when $x \in ( U(\partial K) \backslash K ) \cap C$, it follows that $x_2 > 0$ which can be chosen sufficiently small so that 
\\ ~~ \\
$(x-y)^\top F(x) = (x_1 - y_1) + [0,1] x_2^2 \log (|x_2|)$ 
 \\ 
$ ~~~~~~~~~~~~~ \subset  (x_1 - y_1) + [-1,1] |  x_2^2 \log (|x_2|)| $ 
 \\ 
$~~~~~~~~~~~~~ = (x_1 - y_1) + [-1,1] (x_2 - y_2)^2 | \log (|x_2 - y_2|)|$ \\ 
$ ~~~~~~~~~~~~~ \subset (x - y)^\top F(y) + [-1,1] |x - y|^2 |\log (|x - y|)|$
\\ 
$ ~~~~~~~~~~~~~ = (x - y)^\top F(y) + [-1,1] |x - y| \rho (|x - y|)$, \\~~\\
where $\rho(w):= w \log w$, for all $\omega \geq 0$. The function $\rho$ is Osgood \cite{redheffer1972theorems}; hence, a uniqueness function. 
\end{example}

\subsection{Pre-Invariance Using Locally Lipschitz Barrier Functions} \label{Subsec.2}

Another approach to conclude forward pre-invariance using cone conditions without restricting the regularity of the flow map, consists of replacing the flow condition \eqref{eq.1} in Theorem \ref{prop1} by a cone condition to be satisfied on the external part of a sufficiently small neighborhood of the set $\partial K_e \cap C$. Inspired by \cite[Theorem~5.2.1]{Aubin:1991:VT:120830}, the flow condition in the following statement uses the external contingent cone $E_K$ defined in \eqref{eq.exter}. 

\begin{theorem} \label{prop6}
Given a hybrid system $\mathcal{H} = (C,F,D,G)$ and a barrier function candidate $B$ defining the set $K$ in \eqref{eqK}. The set $K$ is forward pre-invariant if \eqref{eq.2}-\eqref{eq.2rj} hold, and there exists a neighborhood $U(\partial K)$ such that
\begin{equation}
\label{eq.econe}
\begin{aligned} 
\eta \in E_K(x) &~~~~~~~\forall \eta \in F(x) \cap T_C(x)~\mbox{and} \\ & ~~~~~~~\forall x \in ( U(\partial K) \backslash K ) \cap C,
\end{aligned}
\end{equation}
where $E_K$ is introduced in \eqref{eq.exter}.
\end{theorem}

\begin{proof}
The proof that the solutions, under \eqref{eq.2}-\eqref{eq.2rj}, cannot leave the set $K$ after a jump according to scenario \ref{item:S1} is the same as in the proof of Proposition \ref{prop1}. Next, as in \cite[Theorem 5.2.1]{Aubin:1991:VT:120830}, we show that the trajectories starting from the set $\partial K \cap C$ cannot leave the set $K$ by flowing according to scenario \ref{item:S2}. Indeed, assume that a solution $x$ starting from $\partial K \cap C$ leaves the set $K$ by flowing according to scenario \ref{item:S2}. Then, there exists $j \in \mathbb{N}$, such that, for $t'_2$ small enough and for all $t \in (t'_1,t'_2]$, $x(t , j) \in (U(\partial K) \backslash K) \cap C $, where $t'_1$ and $t'_2$ are as in \ref{item:S2}. Furthermore, using 
\ifitsdraft  Lemma \ref{lem1}\else \cite[Lemma 2]{draftautomatica}\fi, we conclude that $ \dot{x}(t,j) \in T_C(x(t,j)) $ for almost all $t \in [t'_1, t'_2]$. Next, since the distance function with respect to the set $K$ is locally Lipschitz and the solution $x(\cdot,j)$ is absolutely continuous on the interval $[t'_1,t'_2]$, it follows that $\delta_K(\cdot) := |(x(\cdot,j)|_K$ is also absolutely continuous on that same interval. Hence, for almost all $t \in [t'_1,t'_2]$, the time derivative $\dot{\delta}_K (t)$ exists and satisfies  
\begin{align*}
\lim_{h \rightarrow 0^+} & \frac{\delta_K(t+h) - \delta_K(t)}{h} = \nonumber \\ &
\liminf_{h \rightarrow 0^+} \frac{|x(t,j) + h \dot{x}(t,j)|_K - |x(t,j)|_K}{h}.  
\end{align*} 
Since $\dot{x}(t,j) \in F(x(t,j)) \cap T_C(x(t,j))$ for almost all $t \in [t'_1, t'_2]$ and 
$ x((t'_1,t'_2],j) \subset (U(\partial K) \backslash K) \cap C $, using \eqref{eq.econe}, we conclude that  $\dot{x}(t,j) \in E_K(x(t,j))$ for almost all $t \in (t'_1,t'_2)$, which implies that $\dot{\delta}_K(t) \leq 0$ for almost all $t \in (t'_1,t'_2)$. Thus, $|x(t'_2,j)|_K - |x(t'_1,j)|_K \leq 0$. The contradiction follows since $x(t'_2,j) \in C \backslash K$ and $x(t'_1,j) \in K$ which means that 
$|x(t'_2,j)|_K - |x(t'_1,j)|_K > 0$.
\ifitsdraft 
\else
\hfill $\square$
\fi
\end{proof}

\begin{remark} \label{remdisc}
The flow condition \eqref{eq.econe} in Theorem \ref{prop6} is similar to the flow condition \eqref{eq.3bis} in Remark \ref{cor1nw}. Indeed, in \eqref{eq.econe}, we are using the distance function $B(x) := |x|_K$ as a barrier function candidate defining the set $K$ according to Definition \ref{def.bar}. However, since the distance function to a set is only locally Lipschitz, the gradient-based inequality in \eqref{eq.3bis} is replaced by the limit in \eqref{eq.exter} that has the same implication on the monotonic behavior of $t \mapsto B(x(t,0))$ on the neighborhood $(U(\partial K) \backslash K) \cap C$. Furthermore, according to \cite[Corollary 5.2.3]{Aubin:1991:VT:120830}, the external cone condition \eqref{eq.econe} in Theorem \ref{prop6} is satisfied provided that, for any $x \in ( U(\partial K) \backslash K ) \cap C$, $F(x) \cap T_C(x) \subset F(y) \subset T_K(y)$, where $y$ is the projection of $x$ on the set $K$. 
\end{remark}

Inspired by the discussion in Remark \ref{remdisc}, in the following statement, we replace the flow condition \eqref{eq.econe} in Theorem \ref{prop6} by a condition involving a general locally Lipschitz barrier function candidate instead of only the distance function. The proof is in     \cite{draftautomatica}. 

\begin{theorem} \label{prop6bis}
Given a hybrid system $\mathcal{H} = (C,F,D,G)$ and a scalar  locally Lipschitz barrier function candidate $B$ defining the set $K$ in \eqref{eqK}. The set $K$ is forward pre-invariant if \eqref{eq.2}-\eqref{eq.2rj} hold, and there exists a neighborhood $U(\partial K)$ such that
\begin{equation}
\label{eq.econebis}
\begin{aligned} 
\displaystyle \max_{\zeta \in \partial_C B(x)} \langle \zeta, \eta \rangle \leq 0 & ~~~~~~~ \forall \eta \in F(x) \cap T_C(x)~\mbox{and} \\ & ~~~~~~~\forall x \in ( U(\partial K) \backslash K ) \cap C,
\end{aligned}
\end{equation}
where $\partial_C B$ is the Clarke generalized gradient of $B$, see 
\ifitsdraft Definition \ref{defgen}. \else \cite[Definition 16]{draftautomatica}. \fi
\end{theorem}

\ifitsdraft
\blue{\begin{proof}
The proof that the solutions, under \eqref{eq.2}-\eqref{eq.2rj}, cannot leave the set $K$ after a jump according to scenario \ref{item:S1} is the same as in the proof of Proposition \ref{prop1}. Next, as in the proof of Theorem \ref{prop6}, we show that the trajectories starting from the set $\partial K \cap C$ cannot leave the set $K$ by flowing according to scenario \ref{item:S2}. Indeed, assume that a solution $x$ starting from $\partial K \cap C$ leaves the set $K$ by flowing according to scenario \ref{item:S2}. Then, there exists $j \in \mathbb{N}$, such that, for $t'_2$ small enough and for all $t \in (t'_1,t'_2]$, $x(t, j) \in (U(\partial K) \backslash K) \cap C$, where $t'_1$ and $t'_2$ are as in \ref{item:S2}. Next, since $B$ is locally Lipschitz and the solution $x(\cdot,j)$ is absolutely continuous on the interval $[t'_1,t'_2]$, it follows that $B(x(\cdot,j))$ is also absolutely continuous on that same interval. Hence, the time derivatives $\dot{B}(x(t,j))$ and $\dot{x}(t,j)$ exist for almost all $t \in [t'_1,t'_2]$ and satisfy 
\begin{align} \label{eq.integ2-bis}
\dot{B}(x(t,j)) = & \lim_{h \rightarrow 0^+} \frac{B(x(t+h,j)) - B(x(t,j))}{h} \nonumber 
\\ = & \lim_{h \rightarrow 0^+} \frac{B(x(t,j) + h \dot{x}(t,j)) - B(x(t,j))}{h} \nonumber \\ 
\leq & \max_{\zeta \in \partial_C B(x(t,j))}   \langle \zeta, \dot{x}(t,j) \rangle.
\end{align} 
The latter inequality is true since 
$B$ is locally Lipschitz, see 
\cite[Page 7]{sanfelice2007invariance} and \cite{clarke1990optimization} for more details. Furthermore, using \ifitsdraft  Lemma \ref{lem1}\else \cite[Lemma 2]{draftautomatica}\fi, we conclude that $ \dot{x}(t,j) \in F(x(t,j)) \cap T_C(x(t,j))$ for almost all $t \in [t'_1, t'_2]$. The latter implies, under \eqref{eq.econebis}, that $\dot{B}(x(t,j)) \leq 0$ for almost all $t \in [t'_1,t'_2]$; thus, $B(x(t'_2,j)) - B(x(t'_1,j)) \leq 0$. Hence, the contradiction follows since $x(t'_2,j) \in C \backslash K$, $x(t'_1,j) \in K$, and $B(x(t'_2,j)) - B(x(t'_1,j)) > 0$.   
\ifitsdraft 
\else
\hfill $\square$
\fi
\end{proof}}
\fi

\begin{remark}
In Theorem \ref{prop6bis}, we consider only the case of scalar barrier function candidates defining the set $K$. Conveniently, when the set $K$ is defined via multiple locally Lipschitz candidates according to \eqref{eqK}, we can then use \eqref{scalbar} to construct a scalar barrier function that is locally Lipschitz and at the same time defines the set $K$ according to \eqref{eqK}. Furthermore, if we compare conditions \eqref{eq.1} and \eqref{eq.econebis} while replacing the scalar function $B$ in \eqref{eq.econebis} by the function $\bar{B}$ in \eqref{scalbar}, we notice that the inequality in \eqref{eq.1}, for each $i \in \{1,2,\ldots,m\}$, is checked in more points compared to \eqref{eq.econebis}. Finally, when the set $C$ is open and $K \subset C$, condition \eqref{eq.07prop0bis} covers what is proposed in \cite[Theorem 2]{glotfelter2017nonsmooth} for unconstrained differential inclusions.
\end{remark}

\begin{example} \label{exp1nwbis}
Consider the hybrid system $\mathcal{H} = (C,F,D,G)$ with $C$ and $G$ as in Example \ref{exp1}, and $F$ and $D$ given by
\begin{align*}
F(x) := & \begin{bmatrix}  - x_2^2 x_1 \\ - (|x|^2 - [0,1]) (2-|x|^2) \end{bmatrix}
~~~\forall x \in C, \\
D := & \left\{ x \in \mathbb{R}^2 : x_2 = 0,~|x| \leq 1 \right\}.
\end{align*} 
We employ Theorem \ref{prop6} to verify forward pre-invariance of the closed set $ K := \left\{ x \in C \cup D : |x|^2 - 1 \leq 0 \right\}$. Indeed, this set admits the scalar barrier function candidate $B(x) := x_2 (|x|^2 - 1)$. Furthermore, for any $x \in D \cap K$, $G(x) \subset C \cup D$; hence, \eqref{eq.2rj} is satisfied. Moreover, for every $x \in K \cap D$ and for every $\eta \in G(x)$, there exists $\alpha \in [0,1]$ such that $\eta = [0~~\alpha |x_1|]^\top$ and 
$B(\eta) = \alpha |x_1| (\alpha^2 |x_1|^2 - 1)$; hence, since $x_1 \leq 1$ for all $x \in K \cap D$, \eqref{eq.2} follows. Furthermore, the set $(U(\partial K) \backslash K) \cap C$ can be chosen to be the open set $(U(\partial K) \backslash K) \cap C = \left\{ x \in C : |x_1|<1,~|x| \in (1, 2) \right\} $. Note that, the function $|\cdot|_K$ satisfies $|x|_K = |x| - 1$ for all $x \in (U(\partial K) \backslash K) \cap C$ and is differentiable on the latter set. Hence, $\liminf_{h \rightarrow 0^+} \frac{|x+hv|_K - |x|_K}{h} = \langle [ x_1 ~~ x_2]^\top, v \rangle / |x|$ and $ E_K(x) = \left\{ v \in \mathbb{R}^2 : \langle [ x_1 ~~ x_2]^\top, v \rangle \leq 0 \right\}$ for all $x \in (U(\partial K) \backslash K) \cap C$. The This implies that \eqref{eq.econe} holds since, for every $\eta \in F(x)$ and for every $x \in (U(\partial K) \backslash K) \cap C$, $\langle [ x_1 ~~ x_2]^\top, \eta \rangle = - x_1^2 x_2^2 - x_2 (|x|^2 - \epsilon) (2-|x|^2) \leq 0$ for some $\epsilon \in [0,1]$. 
\end{example}

\subsection{Pre-Invariance Using a Relaxed Flow Condition} \label{Subsec.3}

Inspired by the property of the uniqueness functions in Definition \ref{defuniq}, we are able to relax the sign of the inequality in \eqref{eq.1} provided that a growth condition involving a uniqueness function and the barrier function candidate is satisfied. Such a relaxation follows the lines of what is proposed in \cite{DAI201762, 10.1007/978-3-642-39799-8_17} for hybrid automata, see Remark \ref{rem.compar}.   

\begin{proposition} \label{pro.nuniq}
Given a hybrid system $\mathcal{H} = (C,F,D,G)$ and a $\mathcal{C}^1$ barrier function candidate $B$ defining the set $K$ in \eqref{eqK}. The set $K$ is forward pre-invariant if \eqref{eq.2}-\eqref{eq.2rj} hold and there exists a neighborhood $U(M_i)$ such that 
\begin{align} \label{eq.07prop0}
\langle \nabla B_i(x), \eta \rangle \leq ~ \rho(B_i(x))~~ & \forall x \in (U(M_i) \backslash K_{ei}) \cap C~\mbox{and} \nonumber  \\ & \forall \eta \in F(x) \cap T_C(x), 
\end{align}
where $\rho : \mathbb{R} \rightarrow \mathbb{R}$ is a uniqueness function.
\end{proposition}

\begin{proof}
Under \eqref{eq.2}-\eqref{eq.2rj}, the proof that the solutions starting from the set $K$ cannot jump outside the set $K$ following scenario \ref{item:S1} is the same as in the proof of Theorem \ref{prop1}. The only remaining way to leave the set $K$ is by flowing according to scenario \ref{item:S2}. We conclude in this case, for $t'_2$ small enough, the existence of $k \in \left\{1,\ldots, m \right\}$ such that $B_k(x(t,j)) > 0$ for all $t \in (t'_1,t'_2]$ and $x((t'_1,t'_2], j) \subset (U(\partial K_k) \backslash K_k) \cap C$, where $t'_1$ and $t'_2$ are as in \ref{item:S2}. Furthermore, using \ifitsdraft  Lemma \ref{lem1}\else \cite[Lemma 2]{draftautomatica}\fi, we conclude that $ \dot{x}(t,j) \in T_C(x(t,j)) $ for almost all $t \in [t'_1, t'_2]$. Hence, using \eqref{eq.07prop0}, we conclude that
\begin{align*} 
\frac{d B_k(x(t,j))}{dt} \leq \rho (B_k(x(t,j)))~~~\mbox{for. a. a}~t \in [t'_1, t'_2]
\end{align*}
with $B_k(x(t'_1,j)) = 0$. Since $\rho$ is a uniqueness function, we conclude that $ B_k(x(t,j)) = 0 $  for all $t \in [t'_1, t'_2]$ and the contradiction follows. 
\ifitsdraft 
\else
\hfill $\square$
\fi
\end{proof}

\begin{remark} \label{remns}
As in Theorem \ref{prop6bis}, when the barrier function candidate $B$ is scalar and locally Lipschitz, the statement of Proposition \ref{pro.nuniq} holds true if we replace \eqref{eq.07prop0} therein by
\begin{align} \label{eq.07prop0bis}
\displaystyle \max_{\zeta \in \partial_C B(x)} \langle \zeta, \eta \rangle \leq ~ \rho(B(x))~~ & \forall x \in (U(K) \backslash K) \cap C~\mbox{and} \nonumber  \\ & \forall \eta \in F(x) \cap T_C(x).
\end{align}
\end{remark}

\subsection{From Pre-Invariance to Invariance}

In the following statement, we show when a forward pre-invariant set $K \subset C \cup D$ is forward invariant. The proof can be found in \cite{draftautomatica}.

\begin{proposition} \label{prop2}
A forward pre-invariant set $K \in C \cup D$ is forward invariant if the solutions cannot escape in finite time inside $K \cap C$ and, for any initial condition in the set $(K \cap \partial C) \backslash D$, a nontrivial flow exists.
\end{proposition}

\ifitsdraft 
\begin{proof}
We first recall that having the set $K$ forward pre-invariant implies that all the solutions starting from the set $K$ cannot flow in $C \backslash K$ and can neither jump outside the set $K$. Moreover, since for all $y \in K \backslash C$, $G(y)$ is nonempty, then $G(y) \subset K$, the solutions starting from $K$ cannot die in the set $K \cap D$. Hence, the only way for a solution to be forward noncomplete is either by dying in the set $(K \cap \partial C) \backslash D$, which is avoided by assumption, or by escaping in finite time inside the set $K \cap C$ which is also not possible by hypothesis. Hence, all the maximal solutions starting from the set $K$ remain in $K$ and are defined on an unbounded hybrid time domain.
\end{proof}
\fi

\begin{remark}
One can guarantee that the solutions do not have a finite escape time inside the set $K \cap C$ when, for example, the set $K \cap C$ is compact or when the flow map $F$ is globally bounded in 
$K \cap C$. 
\end{remark}   

\ifitsdraft
\blue{
\begin{example}
 Using Proposition \ref{prop2}, we are able to extend the conclusions in Example \ref{exp1} and conclude forward invariance of the set $K$. Indeed, in Example \ref{exp1}, we showed that the set $K$ is forward pre-invariant and since it is compact, the solutions starting from $K$ cannot blow-up in finite time. Hence, the forward invariance follows if we show that, for every initial condition in the set $ (K \cap \partial C) \backslash D = \left\{ [1~~0]^\top, [-1~~0]^\top \right\}$, a nontrivial flow exits. Indeed, $F(x) = \left\{ 0 \right\}$ for all $x \in (K \cap \partial C) \backslash D$; hence, the system admits nontrivial constant solutions of the form $x(t,0) = x_o$ for all $t \geq 0$ and $x_o \in (K \cap \partial C) \backslash D$. 
\end{example}}
\fi

\begin{example} [Thermostat]
Using Proposition \ref{prop2}, we extend the conclusions in Example \ref{exptermostat} and show that the set $K$ introduced therein is forward invariant. Indeed, we already showed that the set $K$ is forward pre-invariant. Furthermore, since it is compact, then there is not a possibility of a finite-time escape inside 
$K \cap C$. Finally, we note that $(K \cap \partial C) \backslash D = (K \cap C) \backslash D =
 \left\{ 0 \right\} \times (z_{min},z_{max}] \cup \left\{1 \right\} \times [z_{min},z_{max})$ and, by explicitly solving the flow dynamics on $[0,t_1]$, for some $t_1 > 0$, we conclude that, when $x_o \in \left\{ 0 \right\} \times (z_{min},z_{max}] $, there exists a nontrivial flow given by $q(t,0) = 0,~z(t,0) = (z(0,0) - z_o) e^{-t} +  z_o$. Similarly, when $x_o \in \left\{1 \right\} \times [z_{min},z_{max})$, there exists a nontrivial flow given by $ q(t,0) = 1,~z(t,0) = (z(0,0) - z_o - z_\Delta) e^{-t} +  z_o + z_\Delta $.
\end{example}

In the following result, we propose a qualitative condition implying the existence of a nontrivial solution starting from each element in the set $(K \cap \partial C) \backslash D$.

\begin{proposition} \label{cor.non-pre}
A forward pre-invariant set $K \in C \cup D$ is forward invariant if the solutions cannot escape in finite time inside the set $K \cap C$ and there exists a neighborhood $U(x_o)$ such that 
\begin{align} \label{eq3-++n}
F(x) \cap T_{K \cap C}(x) \neq & ~ \emptyset~~~\forall x \in U(x_o) \cap (K \cap \partial C)~\mbox{and} \nonumber \\ &  ~~~~~ \forall x_o \in (K \cap \partial C) \backslash D.  
\end{align}
\end{proposition}

\begin{proof}
To conclude the proof in this case, we propose to show that
\begin{align}
\label{eqtoprove1}  
T_{K\cap C}(x) \cap F(x) \neq & \emptyset ~~~ \forall x \in U(x_o) \cap \partial (K \cap C)~\mbox{and} \nonumber \\ & ~~~~~~~~ \forall x_o \in (K \cap \partial C) \backslash D.
\end{align}
Indeed, using \eqref{eqtoprove1} and \ifitsdraft Proposition \ref{thm2} \else 
\cite[Proposition 8]{draftautomatica} \fi with the set $K$ therein replaced by $K \cap C$, we conclude the existence of a nontrivial flow starting from each $x_o \in (K \cap \partial C) \backslash D$. Hence, the forward invariance of the set $K$ follows using Proposition \ref{prop2}. Now, in order to show \eqref{eqtoprove1}, we distinguish two complementary situations. First, when $x \in  \partial (K \cap C) \cap \partial C = K \cap \partial C$, in this case, \eqref{eqtoprove1} follows from \eqref{eq3-++n}. Second, when $x \in \partial (K \cap C) \cap \mbox{int}(C)$, in this case, since the set $K$ is forward pre-invariant and since there exist always a nontrivial solution flowing from each $x \in \mbox{int}(C)$ under the standing assumptions, we conclude the existence of a nontrivial solution flowing from $x$ and remaining in $K$ for a nontrivial interval of time. Hence, using \ifitsdraft Proposition \ref{thm1} \else 
\cite[Proposition 7]{draftautomatica}\fi, \eqref{eqtoprove1} follows also when $x \in \partial (K \cap C)  \cap \mbox{int}(C)$, which completes the proof.
\ifitsdraft 
\else
\hfill $\square$
\fi
\end{proof}

\ifitsdraft
\begin{example}
We propose to extend forward pre-invariance of the set $K$, in Example \ref{exp1nw}, in order to conclude forward invariance using Proposition \ref{cor.non-pre}. That is, the set $K$ is compact and  
$(K \cap \partial C) \backslash D  = \left\{ (-1,0),(1,0) \right\} $. Furthermore, for any $ x_o \in (K \cap \partial C) \backslash D$,
\begin{align*}
U(x_o) \cap (K \cap \partial C) = \left\{
\begin{matrix} 
[3/4, 1] \times \left\{0\right\} & \mbox{if}~x_o = (1,0) \\
-[1, 3/4] \times \left\{0\right\} & \mbox{if}~x_o = (-1,0) 
\end{matrix} \right.
\end{align*}
Next, for any $x \in U(x_o) \cap (K \cap \partial C)$, there exits $\delta_x \in [0,1/2)$ with $\delta_{x_o} = 0$ such that $|x|^2 = 1- \delta_x$ and $F(x) \in [0~~[\delta_x-1, \delta_x] (3/4-\delta_x) (1+\delta_x)]^\top $, hence, the element $ v(x) := [0~~\delta_x (3/4-\delta_x) (1+\delta_x)]^\top $ belongs to $F(x) \cap T_{C \cap K} (x)$, hence \eqref{eq3-++n} is satisfied.  
\end{example}
\else
\begin{example} [Bouncing ball]  
For the system in Example \ref{expbb}, we will show forward invariance of the set $K$ using Proposition \ref{cor.non-pre}. Indeed, the set $K$ is compact; hence, the solutions starting from $K$ cannot escape in finite time. Furthermore, $(K \cap \partial C) \backslash D  = \left\{ 0 \right\} \times (0,1]$ and, for any $ x_o \in (K \cap \partial C) \backslash D$, we can take 
$ U(x_o) \cap (K \cap \partial C) = \left\{0 \right\} \times [x_o/2, 1]$. Next, for any $x \in U(x_o) \cap (K \cap \partial C)$ with $U(x_o)$ small enough, we notice that $F(x)$ is transversal to 
$\partial C$ and at the same time $\langle \nabla B(x), F(x) \rangle \leq 0$; hence, using \cite[Corollary 3]{draftautomatica}, we conclude that $F(x) \in D_C(x) \cap T_{K_e}(x)$, which implies, using \cite[Lemma 5]{draftautomatica}, 
that $F(x) \in T_{C \cap K}(x)$ and \eqref{eq3-++n} follows.  
\end{example}
\fi

\section{Sufficient Conditions for Pre-Contractivity and Contractivity Using Barrier Functions} \label{Sec.4}

A pre-contractive set is a forward pre-invariant set such that whenever a solution starts from its boundary, it immediately leaves it and evolves towards its interior. The study of contractive sets is very important since several techniques to derive Lyapunov functions are based on the construction of contractive sets, the resulting Lyapunov functions are known as \textit{set-induced Lyapunov functions}, see \cite{Blanchini:1999:SPS:2235754.2236030, fiacchini2015computation, 272351, giesl2015review}.

\subsection{Definitions}

A definition of contractivity for the particular sets named $C-$sets is proposed in \cite{Blanchini:1999:SPS:2235754.2236030} using the \textit{Minkowskii} functional, also named \textit{gauge} function, for both differential and difference equations, see Definitions 3.3 and 3.4 in \cite{Blanchini:1999:SPS:2235754.2236030}, respectively. The latter approach can be extended for general hybrid inclusions to define pre-contractivity and contractivity for $C-$sets. Indeed, we start recalling that a set $K \subset C \cup D$ is said to be a $C-$set if it is compact, convex and includes the origin in its interior. Moreover, the corresponding Minkowskii functional $\Psi_K: \mathbb{R}^n \rightarrow \mathbb{R}_{\geq 0}$ is given by
\begin{align}  \label{eq.Psipre}
\Psi_K(x) := \inf \left\{ \mu \geq 0 :~x \in \mu K \right\}.
\end{align}

\begin{definition}[Pre-contractivity for $C-$sets] \label{def.cocopre} 
A $C$-set $K \subset C \cup D$ is said to be pre-contractive if 
\begin{align}
\limsup_{h \rightarrow 0^+} & \frac{ \Psi_K(x+ \eta h)-1 }{h} < 0 
\nonumber \\ &
\qquad \forall x \in \partial K \cap C, \quad \forall \eta \in F(x) \cap T_C(x), 
\label{eq.mink1pre}   \\ 
\Psi_K(\eta) & < 1 \qquad \forall x \in D \cap K, \quad \forall \eta \in G(x).  \label{eq.mink2pre} 
\end{align}
\end{definition}

\begin{definition}[Contractivity for $C-$sets] \label{def.cocopre-} 
A $C$-set $K \subset C \cup D$ is said to be contractive if it is pre-contractive and, in addition, starting from each element in the set $(K \cap \partial C) \backslash D$, a nontrivial solution exists.
\end{definition}

The following lemma establishes important consequences of the contractivity properties, in Definitions \ref{def.cocopre} and \ref{def.cocopre-}, on the behavior of the system's solutions. Based on these consequences, we will define pre-contractivity and contractivity for general closed sets that are not necessarily $C-$sets. The proof is in \cite{draftautomatica}.  

\begin{lemma} \label{lemmmm1}
If a $C-$set $K \subset (C \cup D)$ is pre-contractive (respectively, contractive) according to Definition \ref{def.cocopre} (respectively, Definition \ref{def.cocopre-}), then it is forward pre-invariant (respectively, forward invariant) and, for any $x_o \in \partial K$ and any nontrivial solution $x$ starting from $x_o$, there exists $T > 0$ and $J \in \mathbb{N}^*$ such that $x(t,j) \in \mbox{int}(K)$ for all $(t,j) \in \dom x \cap \left[ \left( [0,T] \times \left\{ 0 \right\} \right)  \cup  \left( \left\{ 0 \right\} \times \left\{ 0,1, \ldots, J \right\} \right) \right]$, $ (t,j) \neq (0,0)$.
\end{lemma}

\ifitsdraft
\blue{
\begin{proof}
Given a pre-contractive $C$-set $K \subset C \cup D$, we first establish its forward pre-invariance. Indeed, under \eqref{eq.mink2pre}, all the solutions staring from the set $K$ cannot jump outside the set $K$. Next, we will show that all the solutions starting from  $\partial K$ cannot flow in $C \backslash \mbox{int}(K)$. To this end, we show, using \eqref{eq.mink1pre}, that
\begin{align} \label{eqpr}
\forall x \in \partial K \cap C, ~~~ F(x) \cap T_{C \backslash \mbox{int}(K)}(x) = \emptyset. 
\end{align}
Indeed, this latter fact combined with \ifitsdraft Proposition \ref{thm1} \else 
\cite[Proposition 7]{draftautomatica} \fi implies that there are no solutions starting from 
$\partial K \cap C$ that flow into $C \backslash \mbox{int}(K)$. Thus, forward pre-invariance of the set $K$ follows. To prove \eqref{eqpr}, we use \eqref{eq.mink1pre} to conclude that for each $x \in \partial K \cap C$ and each $\eta \in F(x) \cap T_C(x)$ there exists a sequence $h_k \rightarrow 0$, $k \in \mathbb{N}$, such that $\Psi_K(x + \eta h_k)-1 < 0$ for all $k \in \mathbb{N}$. This property implies that $x + \eta h_k \in \mbox{int}(K)$. Hence, there exists $\epsilon > 0$ such that $x + h_0 (\eta + \epsilon \mbox{int}(\mathbb{B})) \subset K$. Furthermore, using the convexity of the set $K$, we conclude that $x + h (\eta + \epsilon \mbox{int}(\mathbb{B}) ) \subset K$ for all $h \in (0, h_0]$ which implies that $\eta \in D_{\mbox{int}(K)}(x)$ for all $\eta \in F(x) \cap T_C(x)$. Hence, for each $x \in \partial K \cap C$ and $\eta \in F(x)$, either $\eta \notin T_C(x)$ thus 
$\eta \notin T_{C \backslash (\mbox{int}(K) \cap C)}(x) = T_{C \backslash \mbox{int}(K)}(x)$, or $ \eta \in D_{\mbox{int}(K)}(x) $ which implies that $\eta \notin T_{\mathbb{R}^n \backslash \mbox{int}(K)}(x)$, thus, $ \eta \notin T_{C \backslash \mbox{int}(K)}(x) $, which concludes \eqref{eqpr}. Hence, the $C$-set $K$ is forward pre-invariant. 

Next, we show that for any nontrivial solution $x$ starting from $x_o \in \partial K$, there exists $T > 0$ and $J \in \mathbb{N}^*$ such that $x(t,j) \in \mbox{int}(K)$ for all $(t,j) \in \dom x \cap \left[ \left( [0,T] \times \left\{ 0 \right\} \right)  \cup  \left( \left\{ 0 \right\} \times \left\{ 0,1, \ldots, J \right\} \right) \right]$, $ (t,j) \neq (0,0)$. Indeed, under \eqref{eq.mink2pre}, all the possible jumps from $\partial K$ satisfy $x(t,j) \in \mbox{int}(K)$ for all $(t,j) \in (\dom x \cap \left( \left\{ 0 \right\} \times \left\{ 0,1 \right\} \right)) \backslash (0,0)$. Moreover, since we have already shown, under \eqref{eq.mink1pre}, the non existence of  flows in $C \backslash \mbox{int}(K)$, it follows that all the possible flows from $\partial K$ satisfy $x(t,0) \in \mbox{int}(K)$ for all $(t,0) \in \dom x \cap \left( [0,T] \times \left\{ 0 \right\} \right) \backslash \left\{(0,0)\right\}$ for some $T>0$, which concludes the proof when the $C$-set $K$ is pre-contractive. Now, given a contractive $C$-set $K \subset C \cup D$, it is already pre-contractive, hence, forward pre-invariant, and, for each solution $x$ starting from $x_o \in \partial K$, 
if $\dom x \backslash \left\{(0,0)\right\} \neq \emptyset$, then there exists $T>0$ and $J \in \mathbb{N}$ such that, for all $(t,j) \in \dom x \cap \left( [0,T] \times \left\{ 0,1,..., J \right\} \right) \backslash 
\left\{ (0,0) \right\}$, $x(t,j) \in \mbox{int}(K)$. 

It remains only to show that the set $K$ is forward invariant. Since the $C-$set $K$ is compact, there is not a possibility of finite-time escape inside $K$; thus, using Proposition \ref{prop2}, it is enough to show that all the maximal solutions starting from the set $K$ have a nontrivial hybrid time domain. To this end, we distinguish three complementary situations. First, when solutions start from $K \cap D$, we use the fact that $G(x)$ is nonempty for all $x \in D$ to conclude that there is always a possibility of jumps from $K \cap D$. Second, when the solutions start from $ (K \backslash D) \cap \mbox{int}(C)$ there is always a possibility of flowing in $C$ under the standing assumptions; hence, the maximal solutions starting from $ (K \backslash D) \cap \mbox{int}(C)$ cannot be trivial. Third, when a solution starts from $x_o \in (K \backslash D) \cap \partial C$, in this case, using Definition \ref{def.cocopre-} and the fact that $(K \backslash D) \cap \partial C = (K \cap \partial C) \backslash D$, we conclude that the maximal solutions starting from 
$(K \backslash D) \cap \partial C$ are nontrivial. Consequently, all the maximal solutions starting from the set $K$ have a nontrivial hybrid time domain, which concludes the forward invariance of the set $K$.
\ifitsdraft 
\else
\hfill $\square$
\fi
\end{proof}}
\fi

For general closed sets, we cannot use the Minkowskii functional to define the contractivity notions since $K$ may not be convex. Consequently, a trajectory-based definition, in the case of differential inclusions, is proposed in \cite{Aubin:1991:VT:120830} under the name of \textit{strict invariance}. In this section, we propose definitions of contractivity and pre-contractivity for hybrid systems  based on the behavior of the solutions after reaching the boundary of the considered set. The aim of the proposed definitions is to preserve the properties established in Lemma \ref{lemmmm1}. Furthermore, sufficient conditions in terms of barrier function candidates defining the (closed) set $K$ are proposed.

\begin{definition}[Pre-contractivity for general sets]  A closed set $K \subset C \cup D$ is said to be pre-contractive if it is forward pre-invariant and for every $x_o \in \partial K$ and every nontrivial solution $x$ starting from $x_o$, there exists $T > 0$ and $J \in \mathbb{N}^*$ such that $x(t,j) \in \mbox{int}(K)$ for all $(t,j) \in \dom x \cap \left[ \left( [0,T] \times \left\{ 0 \right\} \right)  \cup  \left( \left\{ 0 \right\} \times \left\{ 0,1, \ldots, J \right\} \right) \right]$, $ (t,j) \neq (0,0)$.
\end{definition} 

\begin{definition}[Contractivity for general sets]  A closed set $K \subset C \cup D$ is said to be contractive if it is pre-contractive and forward invariant. 
\end{definition} 

\begin{remark} \label{rem.precont}
It is useful to notice that, in the particular case of differential inclusions, the pre-contractivity of a closed set $K \subset \mathbb{R}^n$ reduces to the nonexistence of a (non-hybrid) solution $t \mapsto x(t)$ starting from any $x_o \in \partial K$ such that 
$x([0,T]) \subset \mathbb{R}^n \backslash \mbox{int}(K)$ for some $T > 0$. 
\end{remark}

\subsection{Pre-Contractivity}

Next, we propose to characterize contractivity notions using barrier functions defining general closed sets\ifitsdraft \else\footnote{The characterization for the case of $C$-sets is treated in the report version \cite{draftautomatica}.}\fi. Our approach is mainly based on \ifitsdraft Lemma \ref{lem2bis} \else \cite[Lemma 4]{draftautomatica}\fi, which characterizes the Dubovitsky-Miliutin cone $D_{\mbox{int}(K)}$ at the boundary of the considered closed set in terms of the barrier function candidate defining the set. Furthermore, the latter fact is combined with \ifitsdraft Theorem \ref{thm3} 
\else 
\cite[Theorem 6]{draftautomatica} 
\fi 
in order to conclude contractivity. 

\begin{theorem} \label{prop5}
Given a hybrid system $\mathcal{H} = (C,F,D,G)$ and a $\mathcal{C}^1$ barrier function candidate $B$ defining the set $K$ in \eqref{eqK}. The set $K$ is pre-contractive if, for each $i = \left\{1,2,\ldots,m \right\}$,
\begin{equation}
\label{eq.7prop}
\begin{aligned}
\langle \nabla B_i(x), \eta \rangle < &~0~~\forall x \in M_i \cap C \\ & ~~~~~~~~~~~\forall \eta \in F(x) \cap T_C(x), 
\end{aligned}
\end{equation}
\begin{align}
F(x) \cap T_{\partial C \cap \partial K}(x) = & ~\emptyset~~\forall x \in \partial K \cap \partial C, \label{eq.7propbis} 
\end{align}
\begin{align} 
B(\eta)< & ~ 0~~~~~~~~~~~~\forall x \in  K \cap D~~~\forall  \eta \in G(x), \label{eq.8prop} \\
G(x) \subset & ~C \cup D~~~~~~~~~~~~\forall x \in K \cap D,  \label{eq.8proprj-} \\ 
G(x) \subset & ~\mbox{int}(C \cup D)~~~~~~\forall x \in \partial K \cap D. \label{eq.8proprj}
\end{align}  
\end{theorem}

\begin{proof}
We consider, without loss of generality, the hybrid extension of $\mathcal{H}$ denoted by $\mathcal{H}_e := (C_e,F_e,D,G)$ where $C_e := U(C) \subset \mathbb{R}^n$ and $ F_e : C_e \rightrightarrows \mathbb{R}^n$ is any extension of $F$ to $C_e$ that preserves the standing assumptions. Also, we recall that $K_e = \left\{ x \in \mathbb{R}^n : B(x) \leq 0 \right\}$ and $K_{ei} := \left\{ x \in \mathbb{R}^n : B_i(x) \leq 0 \right\}$; hence, $K_e = \cap^{m}_{i=1} K_{ei}$. We start noticing that for each $x_o \in \partial K_e \cap C$ there exists $I_{x_o} \subset \left\{1,2,...,m \right\}$ such that $x_o \in M_i$ if and only if $i \in I_{x_o}$. That is, we have $x_o \in \cap_{i \in I_{x_o}} (\partial K_{ei} \cap C)$ and, also, $x_o \in \cap_{i \in I_{x_o}} (\partial K_{ei} \cap \mbox{int}(C_e) )$. Next, using \ifitsdraft Lemma \ref{lem2bis} \else \cite[Lemma 4]{draftautomatica} \fi under \eqref{eq.7prop}, it follows that $ F_e(x_o) \cap T_C(x_o) \subset D_{\mbox{int}(K_{ei})}(x_o)$ for all $i \in I_{x_o}$. Moreover, when $i \notin I_{x_o}$, we also have $F_e(x_o) \cap T_C(x_o) \subset D_{\mbox{int}(K_{ei})}(x_o)$ since in this case $ D_{\mbox{int}(K_{ei})}(x_o) = \mathbb{R}^n$. Hence, we conclude that $F_e(x_o) \cap T_C(x_o) \subset D_{\mbox{int}(K_{e})}(x_o)$ since $\cap^{m}_{i=1} D_{\mbox{int}(K_{ei})}(x_o) = D_{\mbox{int}(\cap^{m}_{i=1} K_{ei})}(x_o)$. 

As a second step, we show that $F(x_o) \cap T_{C \backslash \mbox{int}(K_e)}(x_o) = \emptyset$. Indeed, for each $x \in \partial K_e \cap C$ and $\eta \in F(x)$, either $\eta \notin T_C(x)$, hence $\eta \notin T_{C \backslash (\mbox{int}(K_e) \cap C)}(x) = T_{C \backslash \mbox{int}(K_e)}(x)$, or $ \eta \in D_{\mbox{int}(K_e)}(x) $ which implies that $\eta \notin T_{\mathbb{R}^n \backslash \mbox{int}(K_e)}(x)$ thus $ \eta \notin T_{C \backslash \mbox{int}(K_e)}(x) $. Now, using \ifitsdraft Proposition \ref{thm1} \else \cite[Proposition 7]{draftautomatica}\fi, we conclude the non existence of any solution $x$ to $\mathcal{H}_e$ starting from $x_o$ and flowing in $C \backslash \mbox{int}(K_e)$ along a nontrivial time interval. Hence, if a nontrivial flow $x$ exists starting from $x_o$, then $x(t,0) \in C_e \backslash (C \backslash \mbox{int}(K_e))$ for all $(t,0) \in \dom x \cap \left([0,T] \times 0 \right) \backslash (0,0)$ and for some $T>0$. That is, for each solution $x$ to $\mathcal{H}_e$ flowing from $x_o$, either there exists $T>0$ such that $x((0,T],0) \subset \mbox{int}(K_e)$ or, there exists $T>0$ such that $x((0,T],0) \subset C_e \backslash C$. Particularly, when $x_o \in \partial K_e \cap \mbox{int}(C) = \partial K \cap \mbox{int}(C) $, a nontrivial flow always exists, hence, $x((0,T],0) \in \mbox{int}(K)$ for some $T>0$. Furthermore, when $x_o \in  \partial C \cap \partial K$, using \ifitsdraft Proposition \ref{thm1} \else \cite[Proposition 7]{draftautomatica} \fi under \eqref{eq.7propbis}, we conclude that $x((0,T],0) \cap (\partial K \cap \partial C) = \emptyset$ for some $T>0$. In other words, the set $\partial K \cap \partial C$ is not weakly forward invariant under \eqref{eq.7propbis}.  Hence, using the fact that $\partial K \cap C = ( \partial K_e \cap C ) \cup (\partial K \cap \partial C) $, it follows that each solution $x$ flowing from $x_o \in \partial K \cap C$  there exits $T>0$ such that either $x((0,T],0) \subset \mbox{int}(K_e) \backslash ( \partial C \cap \partial K ) \cap C$, or $x((0,T],0) \subset C_e \backslash C$. Going back to $\mathcal{H}$, the latter scenario is excluded; hence, for any solution flowing from $x_o \in \partial K \cap C$, there exits $T>0$ such that, if $\dom x \backslash (\mathbb{R}_{\geq 0},0) \neq \emptyset$, then $x(t,0) \subset \mbox{int}(K)$ for all $(t,0) \in \dom x \cap \left( [0,T] \times \left\{ 0 \right\} \right) \backslash (0,0)$. On the other hand, under \eqref{eq.8prop}, \eqref{eq.8proprj-}, and the continuity of $B$, we conclude that  for any 
$x_o \in K \cap D$, $G(x_o) \subset \mbox{int}(K_e) \cap (C \cup D) \subset K$; hence, all the possible jumps from $K \cap D$ maintain the solution in the set $K$. Furthermore, using \eqref{eq.8proprj}, we conclude that, for each $x_o \in \partial K \cap D$ and after any possible jump, the solutions jump to the interior of the set $K$. Hence, if a jump is possible when starting from $x_o \in \partial K \cap D$, then $x(t,j) \in \mbox{int}(K)$ for all $(t,j) \in \dom x \cap \left(\left\{ 0 \right\} \times \left\{0,1\right\} \right) \backslash (0,0) \neq \emptyset$; which completes the proof. 
\ifitsdraft 
\else
\hfill $\square$
\fi
\end{proof}

\begin{example} \label{exprj}
Consider the hybrid system $\mathcal{H}$ given by 
\\
$F(x) := \begin{bmatrix} - (x_2 + 1) \\ - 2 (x_2 + 1) +  x_1  \end{bmatrix}~~~\forall x \in C$, 
\\
$ C := \left\{ x \in \mathbb{R}^2 : x_2 \in [0,1],~
|x_1| \leq \sqrt{3}  \right\}$, 
\\
$G(x) := \frac{1}{\sqrt{3}} \begin{bmatrix} x_1 \\ 
\frac{\sqrt{3}}{2}  \end{bmatrix}~~~\forall x \in D$,\\ 
$ D := \left\{ x \in \mathbb{R}^2 : x_2 = 0,~|x_1| \leq \sqrt{3} \right\}$. \\ ~~ \\
We will show that the set 
$$ K := \left\{ x \in \mathbb{R}^2 : x_1^2 + (x_2+1)^2 \leq 4,~~x_2 \geq 0 \right\}, $$ 
which is not a $C$-set, is pre-contractive. Indeed, the set $K$ can be defined using the $\mathcal{C}^1$ barrier function candidate
$ B(x) := [x_1^2 + (x_2+1)^2 - 4~~-x_2]^\top $. 
To verify \eqref{eq.8prop}-\eqref{eq.8proprj}, we note that, for each $x \in K \cap D$, $G(x) = [ \frac{x_1}{\sqrt{3}} ~~ \frac{1}{2} ]^\top \in \mbox{int}(C \cup D) $ since $(1 / \sqrt{3}) x_1 \leq 1$ for all $x \in K \cap D = \left\{ x \in \mathbb{R}^2 : x_2 = 0,~~|x_1| \leq \sqrt{3} \right\}$; hence  \eqref{eq.8proprj-}-\eqref{eq.8proprj} hold. Moreover,  $B(G(x)) = [ x_1^2/3 - 7/4 ~~ -1/2 ]^\top < 0$ since $|x_1| \leq \sqrt{3}$ for all $x \in K \cap D$; hence, \eqref{eq.8prop} also follows.  Furthermore, $\langle \nabla B_1(x), F(x) \rangle  = - 4(x_2 + 1)^2 < 0 $ for all $ x \in C \cap M_1 = \left\{ x \in \mathbb{R}^2 : x_2 \geq 0,~|[x_1~~x_2+1]| = 2 \right\} $. 
Next, for every  
$ x \in M_2 \cap C = \left\{ x \in \mathbb{R}^2 : x_2 = 0,~|x_1| \leq \sqrt{3} \right\}$, 
$F(x) =  [ -1 ~~ - 2 + x_1 ]^\top \notin T_C(x)$ since $ F_2(x)  = -2 + x_1 < 0 $ for all $x \in M_2 \cap C$; hence, \eqref{eq.7prop} is satisfied. To show that 
\eqref{eq.7propbis} holds, we note that 
$\partial K \cap \partial C = \left\{ [0~~1]^\top \right\} \cup M_2$. Furthermore, when $x = [0~~1]^\top$, $ F(x) = [-2~~-4]^\top \notin T_{\partial K \cap \partial C}(x) = 0$, and, when $x \in M_2$, we have already shown that $F(x) \notin T_C(x)$. The two latter facts imply that $F(x) \notin T_{\partial C \cap \partial K}(x)$ for all $x \in \partial K \cap \partial C$. Hence, 
pre-contractivity of the set $K$ follows according to Theorem \ref{prop5}. 
\end{example}

\ifitsdraft
When the set $K$ is defined using a scalar barrier function candidate, Theorem \ref{prop5} reduces to the following statement.
\begin{corollary} \label{prop5sc}
Consider a $\mathcal{C}^1$ scalar barrier function candidate defining the closed set $K$ as in \eqref{eqK}. The set $K$ is pre-contractive if \eqref{eq.7propbis}-\eqref{eq.8proprj} hold and 
\begin{align}
\langle \nabla B(x), \eta \rangle < &~0~~\forall x \in \partial K_e \cap C~~\forall \eta \in F(x) \cap T_C(x), \label{eq.7propsc} 
\end{align}
\end{corollary}	
\fi

The following result extends Theorem \ref{prop5} when the barrier function candidate defining the closed set $K$ is only locally Lipschitz. The proof is in \cite{draftautomatica}. 

\begin{theorem} \label{prop5lip}
Given a hybrid system $\mathcal{H} = (C,F,D,G)$ and a scalar locally Lipschitz barrier function candidate $B$ defining the set $K$ in \eqref{eqK}. The set $K$ is pre-contractive if \eqref{eq.7propbis}-\eqref{eq.8proprj} hold and
\begin{align}\label{eq.7proplip} 
\hspace{-0.4cm} \displaystyle \max_{\zeta \in \partial_C B(x)} \langle \zeta, \eta \rangle < 0~~\forall x \in K_e \cap C,~
\forall \eta \in F(x) \cap T_C(x).
\end{align}
\end{theorem}

\ifitsdraft
\blue{
\begin{proof}
Using Theorem \ref{prop6bis} under \eqref{eq.8prop}-\eqref{eq.8proprj-} and 
\eqref{eq.7proplip}, we conclude that the set $K$ is forward 
pre-invariant. Furthermore, the proof that, under \eqref{eq.8prop} and \eqref{eq.8proprj}, the solutions jumping from $\partial K \cap D$, jump to $\mbox{int}(K)$ and that, under \eqref{eq.7propbis}, the solutions cannot flow in $\partial K \cap \partial C$ for a nontrivial interval of time is the same as in the proof of Theorem \ref{prop5}. Hence, it remains to show, under \eqref{eq.7proplip}, that every nontrivial solution, flowing from $K_e \cap C$, flows immediately to the interior of the set $K_e$. Indeed, assume that the opposite is true; namely, there exist $t'_2 > t'_1 \geq 0$, $j \in \mathbb{N}$, and a solution $x$ such that $([t'_1,t'_2] \times \left\{j \right\}) \subset \dom x $ and $x( [t'_1,t'_2], j) \subset \partial K_e$. Next, as in the proof of Theorem \ref{prop6bis}, since $B$ is locally Lipschitz and the solution $x(\cdot,j)$ is absolutely continuous on the interval $[t'_1,t'_2]$, it follows that $B(x(\cdot,j))$ is also absolutely continuous on that same interval. Hence, the time derivative $\dot{B}(x(t,j))$ and 
$\dot{x}(t,j)$ exist for almost all $t \in [t'_1,t'_2]$ and satisfy \eqref{eq.integ2-bis}. Furthermore, using \ifitsdraft  Lemma \ref{lem1} \else \cite[Lemma 2]{draftautomatica}\fi, we conclude that 
$\dot{x}(t,j) \in F(x(t,j)) \cap T_C(x(t,j))$ for almost all $t \in [t'_1, t'_2]$. The latter implies, under \eqref{eq.7proplip}, that $\dot{B}(x(t,j)) < 0$ for almost all $t \in [t'_1,t'_2]$; thus, $B(x(t'_2,j)) - B(x(t'_1,j)) < 0$. Contradiction follows since $x(t'_2,j) \in K_e$ and $x(t'_1,j) \in K_e$ yields 
$B(x(t'_2,j)) - B(x(t'_1,j)) = 0$. 
\ifitsdraft 
\else
\hfill $\square$
\fi
\end{proof}}
\fi

\subsection{From Pre-Contractivity to Contractivity}

In the sequel, we complement the sufficient conditions in Theorem \ref{prop5} and Theorem \ref{prop5lip} to conclude contractivity rather than only pre-contractivity.

\begin{proposition} \label{prop5cont}
A pre-contractive closed set $K \subset C \cup D$ is contractive provided that its solutions from $K$ do not escape in finite time inside $K \cap C$ and, starting from each initial condition in $(K \cap \partial C) \backslash D$, a nontrivial flow exists.
\end{proposition}

\ifitsdraft
\begin{proof}
Having the set $K$ pre-contractive, thus forward pre-invariant, we use Proposition \ref{prop2} to conclude that, in the absence of finite-time blow up inside $K \cap C$ and under the existence of a nontrivial flow starting from each initial condition in $ ( K \cap \partial C ) \backslash D$, the set $K$ is forward invariant. The contractivity is proved, since the set $K$ is already 
pre-contractivity.
\ifitsdraft 
\else
\hfill $\square$
\fi
\end{proof}
\fi

In the following statement, the existence of a nontrivial flow starting from $ (K \cap \partial C) \backslash D$ is guaranteed provided that a qualitative tangentiality condition holds. 

\begin{proposition} \label{remrjcor}
Consider a $\mathcal{C}^1$ barrier function candidate defining the set $K$ as in \eqref{eqK}. The set $K$ is contractive if conditions \eqref{eq.7prop}-\eqref{eq.8proprj} hold and there exists a neighborhood $U(x_o)$ such that  
\begin{equation}
\label{eq.rjrj}
\begin{aligned} 
F(x) \cap T_{C}(x) \neq &~ \emptyset ~~~~ \forall x \in U(x_o) \cap K \cap \partial C, \\ & ~~~~~~~ \forall x_o \in 
(K \cap \partial C) \backslash D.
\end{aligned}
\end{equation}
\end{proposition}

\begin{proof}
To conclude the proof, we show that
\begin{equation}
\label{eqtoprove} 
\begin{aligned} 
T_{K\cap C}(x) \cap F(x) \neq \emptyset & ~~ \forall x \in U(x_o) \cap \partial (K \cap C) \\ & 
~~~~  \forall x_o \in (K \cap \partial C) 
\backslash D.
\end{aligned}
\end{equation}
Using \ifitsdraft Proposition \ref{thm2} \else \cite[Proposition 8]{draftautomatica}\fi, \eqref{eqtoprove} implies the existence of a nontrivial flow starting from each $x_o \in (K \cap \partial C) \backslash D$. Under the stated assumptions, the latter fact allows to conclude contractivity of the set $K$ using Proposition \ref{prop5cont}. Now, to prove \eqref{eqtoprove}, we distinguish three complementary situations. 
First, when $x \in  \partial (K \cap C) \cap \partial C  \cap \mbox{int}(K_e) = \partial C \cap \mbox{int}(K_e) $. In this case, \eqref{eqtoprove} follows from \eqref{eq.rjrj} since in this case $T_C(x) = T_{K \cap C}(x)$. 
Second, when $ x \in  \partial (K \cap C) \cap \partial C \cap \partial K_e = \partial C \cap \partial K_e $, in this case, we use \eqref{eq.rjrj} to conclude the existence of $\eta \in F(x)$ such that $\eta \in T_C(x)$. Furthermore, we use \eqref{eq.7prop} under \ifitsdraft Lemma \ref{lem2bis} \else \cite[Lemma 4]{draftautomatica}\fi, to conclude that $\eta \in \mbox{int} (T_{K_e}(x)) = D_{\mbox{int}(K_e)}(x)$, which implies, using \ifitsdraft Lemma \ref{lem6} \else \cite[Lemma 5]{draftautomatica}\fi, that $ \eta \in T_{K \cap C}(x) $. Finally, consider $x \in \partial (K \cap C) \cap \mbox{int}(C) $. In this case, using Theorem \ref{prop5}, the set $K$ is pre-contractive, hence, forward pre-invariant. Combining the latter fact to the existence of a nontrivial flow starting from each $x \in \mbox{int}(C)$ under our standing assumptions, the existence of a nontrivial flow starting from $x$ and remaining in $K$ follows. Thus, under \ifitsdraft Proposition \ref{thm1} \else \cite[Proposition 7]{draftautomatica}\fi, \eqref{eqtoprove} follows when $x \in \partial (K \cap C)  \cap \mbox{int}(C)$.
\ifitsdraft 
\else
\hfill $\square$
\fi
\end{proof}

\begin{example} \label{exprj-}
Using Proposition \ref{remrjcor}, We will show that the closed set $K$ introduced in Example \ref{exprj} is contractive. To this end, we need to show that there is not a possibility of finite-time escape inside $K$. This is the case since the set $K$ is compact. Finally, to show that \eqref{eq.rjrj} holds for all $x \in (K \cap \partial C) \backslash D = \left\{(0,1)\right\}$, 
we note that we can choose $U(x) \cap K \cap \partial C = \{x:=[0~~1]^\top\}$ on which $F(x) = [-2 ~~ -4 ]^\top \in T_C(x)$ since $F_2(x) = -4 < 0$.
\end{example}	

 \section{Conclusion}
 
 The first part of this paper proposed sufficient conditions for forward invariance of closed sets for hybrid systems modeled as hybrid inclusions. The considered closed sets are defined using barrier function candidates and the proposed sufficient conditions in terms of the latter barrier functions are infinitesimal inequalities; namely, not involving any knowledge about the system's solutions, guaranteeing that the set of points on which all the components of the barrier function candidate are nonpositive is forward invariant. Studying forward invariance in the general context of hybrid inclusions offered many technical challenges that have not been handled in the existing literature, to the best of our knowledge. Those challenges are mainly due to the continuous-time evolution of the hybrid inclusion being not necessarily defined on an open set. Hence, elements around the intersection between the zero-level sets of the barrier candidate and the boundary of the set where the continuous-time evolution is defined needed a particular treatment. Finally, we proposed relaxed flow conditions compared to the existing literature by using the notion of uniqueness functions. In the following, we recap the different sufficient conditions for forward pre-invariance proposed in this paper. Given a hybrid system 
$\mathcal{H} = (C,F,D,G)$ and a barrier function candidate $B$ defining the closed set $K$ in \eqref{eqK}. The set $K$ is forward pre-invariant for 
$\mathcal{H}$ if \eqref{eq.2}-\eqref{eq.2rj} hold and either 
\begin{enumerate}[label={\arabic*.},leftmargin=*]
\item \eqref{eq.econe} holds, or
\item $B \in \mathcal{C}^1$ and \eqref{eq.1} holds, or
\item  $B \in \mathcal{C}^1$ and \eqref{eq.07prop0} holds, or
\item $B$ is locally Lipschitz and \eqref{eq.econebis} holds, or
\item  $B$ is locally Lipschitz and \eqref{eq.07prop0bis} holds, or
\item $B \in \mathcal{C}^1$, Assumption \ref{ass1}, \eqref{eq.07prop}, and \eqref{eq.9} hold, and either \ref{item:a}, \ref{item:b}, or \ref{item:c} in Theorem \ref{propnou} holds. 
\end{enumerate}
In the second part of the paper, following the same approach as in the first part, we proposed sufficient conditions for contractivity, which is a stronger property than forward invariance.
 
In the future, it would be  interesting to analyze the necessity of the different sufficient conditions proposed in this paper or to propose new necessary and sufficient ones. Also, it would be interesting to analyze the notions considered in this paper in the presence of perturbations. Furthermore, this work constitutes an important step to analyze the mixed safety plus convergence problem in hybrid systems. Indeed, the latter problem is solved if we show forward invariance of the safety region plus contractivity of the reachable set from the safety region towards a given target. Investigating sufficient (infinitesimal) conditions to guarantee the latter mixed safety-convergence task in hybrid systems is part of our current research efforts.

\bibliography{biblio}
\bibliographystyle{unsrtnat}

\ifitsdraft
\section*{Appendix}

\subsection{Sufficient conditions for Pre-contractivity for $C-$sets}

In the current section, we propose necessary and sufficient conditions for pre-contractivity in terms of barrier function candidates defining a $C$-set. 
\begin{theorem} \label{prop4pre}
Given a hybrid system 
$\mathcal{H}= (C,F,D,G)$, a $C-$set $K \subset \mbox{int}(C \cup D)$ is pre-contractive if and only if there exists a Lipschitz continuous barrier function candidate $B$ defining the $C$-set $K$ as in \eqref{eqK} such that
\begin{equation}
\label{eq.7iffpre} 
\begin{aligned}
\limsup_{h \rightarrow 0^+} & \frac{ B_i(x + \eta h) }{h} < 0 ~~~ \forall x \in M_i \cap C 
 \\ & 
\forall i \in \left\{1,\ldots,m \right\}  ~~ \forall \eta \in F(x) \cap T_C(x), 
\end{aligned}
\end{equation}
\begin{align}
B(\eta) < & ~ 0 ~~~~~~~~~~ \forall x \in K \cap D ~~~ \forall \eta \in G(x), \label{eq.8iffpre}  \\
G(x) \subset & ~ C \cup D ~~~ \forall x \in K \cap D.   \label{eq.8iffprerj}
\end{align} 
\end{theorem}

\begin{proof}
The necessary part of the proof is rather simpler as it relies on translating the Minkowskii functional 
$\Psi_K(x)$ to obtain a barrier function candidate satisfying \eqref{eq.7iffpre}-\eqref{eq.8iffprerj}. That is, the resulting barrier function candidate is given by 
\begin{align} \label{eqconst}
B(x) := & \Psi_K(x)-1,
\end{align}
which is convex hence Lipschitz continuous.

To prove the sufficient part, we, first, pick $x \in K \cap D$. Having $B(\eta) < 0$ for all $\eta \in G(x)$ implies, under the continuity of $B$ and the fact that $K \subset \mbox{int}(C \cup D)$, that either $\eta \in \mbox{int}(K)$ or $\eta \notin C \cup D$. The latter case is not possible under \eqref{eq.8iffprerj}. Hence,  $\eta \in \mbox{int}(K)$ which implies \eqref{eq.mink2pre} since the $C$-set $K$ is convex.

Finally, we prove that \eqref{eq.7iffpre} implies \eqref{eq.mink1pre} using a contradiction. Assume the existence of $x \in \partial K \cap C$ and $\eta \in F(x) \cap T_C(x)$ such that 
\begin{align} \label{eq.eggg}
\limsup_{h \rightarrow 0^+} \frac{ \Psi_K(x+ \eta h)-1 }{h} \geq 0. 
\end{align} 
On the other hand, let $I_x \subset \left\{1,...,m \right\}$ such that $B_i(x) = 0$ if and only if $i \in I_x$. Since $B$ is continuous and $x \in \partial K \cap C$, then $B_j(x) < 0$ for all $j \notin I_x$. Furthermore, using \eqref{eq.7iffpre}, we conclude the existence of a sequence $h_k > 0$, $k \in \mathbb{N}$ such that $h_k \rightarrow 0$ and $B_i(x+h_k \eta) < 0$ for all $k \in \mathbb{N}$ and $i \in I_x$. Since the set $K$ is convex, we conclude that $B_i(x+h \eta) < 0$ for all $h \in (0, h_0]$ and for all $i \in I_x$. Next using the continuity of $B$, we conclude that for $h_0$ sufficiently small, we also have $B_j(x+h \eta) < 0$ for all $h \in (0,h_0]$ and $j \notin I_x$. Thus, $B(x+h \eta) < 0$ for all $h \in (0,h_0]$, which implies that $x + h \eta \in \mbox{int}(K_e)$ for all $h \in (0,h_0]$. Furthermore, since $K \subset \mbox{int}(C \cup D)$, we conclude that for $h_0$ small enough, $x + h \eta \in \mbox{int}(K_e) \cap \mbox{int}(C \cup D) = \mbox{int}(K)$ for all $h \in (0,h_0]$. The latter fact implies that $\Psi_K(x+\eta h) < 1$ for all $h \in (0, h_0]$. Using \eqref{eq.eggg}, we conclude that
\begin{align}\label{eq.eg}
\limsup_{h \rightarrow 0^+} \frac{ \Psi_K(x+ \eta h)-1 }{h} = 0. 
\end{align}  
Furthermore, since the set $K$ is convex and using the first order homogeneity of the Minkowski functional, we conclude the convexity of the function $ \Psi^o_K(h) := \Psi_K(x+h \eta)$ on the interval $h \in [0,h_0]$. That is, 
\begin{align} 
& \limsup_{h \rightarrow 0^+} \frac{ \Psi_K(x+ \eta h)-1 }{h} =  \limsup_{\lambda \rightarrow 0^+} \frac{ \Psi^o_K( h_0 \lambda)-1 }{h_0 \lambda} \nonumber \\ \leq & 
 \limsup_{\lambda \rightarrow 0^+} \frac{ \lambda \Psi^o_K( h_0) + (1-\lambda) -1 }{h_0 \lambda} \leq \frac{ \Psi^o_K( h_0) - 1 }{h_0} < 0
 \end{align}
for all $h \in [0,h_0]$ and for $\lambda := h/h_0$. Hence, the contradiction follows. 
\ifitsdraft 
\else
\hfill $\square$
\fi
\end{proof}

\begin{remark}
The equivalence in the previous statement is shown in the particular case where the $C-$set $K$ satisfies $K \subset \mbox{int}(C \cup D)$. However, the same result remains valid, under the same proof, when $K$ is a $C-$set satisfying $K \subset ( \mbox{int}(C) \cup D ) \backslash (\partial C \cap \partial D)$ and the following extra jump condition holds:
\begin{align}
B(\eta) \nless & 0~~\forall \eta \in G(x) \cap \partial(C \cup D),~~\forall x \in K \cap D.\label{eq.re2}
\end{align}
Indeed, having \eqref{eq.re2} satisfied is important to conclude that, under \eqref{eq.8iffpre}-\eqref{eq.8iffprerj}, the solutions starting from $\partial K$ cannot jump towards $\partial K \cap \mbox{int}(K_e) \subset \partial (C \cup D)$ since it is possible to have $B(x) < 0$ while $x \in \mbox{int}(K_e) \cap \partial (C \cup D) \subset \partial K$.
\end{remark}

\begin{remark}\label{remrj}
In the general case where $K \not\subset ( \mbox{int}(C) \cup D ) \backslash (\partial C \cap \partial D)$, we cannot guarantee for a nontrivial solution flowing from $x_o \in K \cap \partial C$ to satisfy $x([0,\epsilon], x_o) \subset \mbox{int}(K)$, for some 
$\epsilon > 0$, since the solution can flow in $\partial C$ while remaining in $\mbox{int}(K_e)$. In other words, the barrier function candidate does not define the set $K$ on any neighborhood of $C \cup D$ as opposed to the Minkowskii functional which defines the set $K$ in $\mathbb{R}^n$. Therefore, in order to extend Theorem \ref{prop4pre} to the general case where $K \subset C \cup D$, we need to guarantee, additionally, that there is not a possibility of flowing in $\partial C \cap \partial K$ while flowing in $\mbox{int}(K_e)$. The latter fact cannot be characterized in terms of a general barrier function candidate defining the set according to \eqref{eqK}.
\end{remark}

\begin{example} \label{expCsets}
 Consider the hybrid system
\begin{align*}
C := & \left\{ x \in \mathbb{R}^2 : x_2 \geq -1 \right\}, \\
F(x) := & \begin{bmatrix} - [1,2] x_1 + (1/2) x_2 \\ - x_2 - (1/2) x_1 \end{bmatrix}~~~\forall x \in C, \\
D := & \left\{ x \in \mathbb{R}^2 : x_2 \leq -1 \right\} \backslash  \left\{ [-1~~-1]^\top \right\}, \\
G(x) := & [0, 1/2] \begin{bmatrix} x_1 \\ - x_2  \end{bmatrix}~~~\forall x \in D.
\end{align*} 
We would like to study the pre-contractivity of the $C$-set $ K := \left\{ x \in C \cup D : x_1^2 + x_2^2 \leq 2,~x_2 \geq -1 \right\}$ admitting the $\mathcal{C}^1$ multiple barrier function candidate 
$ B(x) := [ |x|^2 - 2 ~~ -(x_2 + 1)]^\top $.  That is, we start noticing that the $C$-set $K$ satisfies $K \subset \mbox{int}(C \cup D)$; hence, Theorem \ref{prop4pre} is applicable. Indeed, since the candidate $B$ is continuously differentiable, we conclude that $\limsup_{h \rightarrow 0^+} \frac{ B_i(x + \eta h) }{h} = \langle \nabla B_i(x), \eta \rangle$ for all $i \in \left\{1,2\right\}$. Furthermore, $ \langle \nabla B_1(x), \eta \rangle \in  [- x_1^2 + x_2^2, -2 x_1^2 + x_2^2] \subset \mathbb{R}_{<0} $ for all $\eta \in F(x)$ and for all $x \in \partial K_1 \cap C = \left\{ x \in \partial K : x_2 \geq -1 \right\}$. Similarly, $ \langle \nabla B_2(x), \eta \rangle = x_2 + (1/2) x_1 = -1 + (1/2) x_1 \leq - 1/2 < 0 $ for all $x \in \partial K_2 \cap C = \left\{ x \in \mathbb{R}^2 : x_2 = -1,~~|x_1| \leq 1 \right\}$ and for all $\eta \in F(x)$. Moreover, $G(x) = [0, 1/2] [ x_1 ~~ 1]^\top \subset [-1/2, 1/2] \times [0, 1/2]$ for all $x \in K \cap D = \left\{ x \in \mathbb{R}^2 : x_2 = -1,~~|x_1| \leq 1 \right\} $; hence, $B_1(\eta) < 0$ and $B_2(\eta) < 0$ for all $\eta \in G(x)$ and for all $x \in K \cap D$. Thus, pre-contractivity of the set $K$ follows using Theorem \ref{prop4pre}. 
\end{example}

\ifitsdraft
A statement similar to Theorem \ref{prop4pre} can be deduced using only a scalar barrier function candidate $B$.
\begin{corollary} \label{prop4pre-}
A $C-$set $K \subset \mbox{int}(C \cup D)$ is pre-contractive if and only if there exists a Lipschitz continuous scalar barrier function candidate $B$ defining the set $K$ as in \eqref{eqK} such that 
\eqref{eq.8iffpre}-\eqref{eq.8iffprerj} are satisfied and 
\begin{align} \label{eq.7iffpre..} 
\limsup_{h \rightarrow 0^+}  \frac{ B(x + \eta h) }{h} < 0 & ~~~~ \forall x \in \partial K \cap C  \nonumber \\ & \qquad ~~\forall \eta \in F(x) \cap T_C(x).
\end{align}
\end{corollary}
\begin{proof}
the sufficient part of the statement is a straightforward consequence Theorem \ref{prop4pre}. Furthermore, the necessary part is also true since the constructed barrier function in \eqref{eqconst} in the proof of Theorem \ref{prop4pre} is a scalar one.
\ifitsdraft 
\else
\hfill $\square$
\fi
\end{proof}
\fi

\subsection{From pre-contractivity to contractivity in the case of $C-$sets}

The previous sufficient conditions can be complemented in order to conclude contractivity rather than only 
pre-contractivity. That is, in the following, we propose sufficient qualitative conditions allowing the existence of nontrivial flows starting from any element in the set $(K \cap \partial C) \backslash D$ as required in Definition \ref{def.cocopre-}.
\begin{proposition} \label{lemconta}
A $C$-set $K \subset C \cup D $ is contractive if it is pre-contractive and
\begin{equation}
\label{eqexists}
\begin{aligned} 
F(x) \cap T_{C}(x) \neq  \emptyset~~~~\forall & x \in U(x_o) \cap K \cap \partial C, \\  
\forall & x_o \in (K \cap \partial C) \backslash D.
\end{aligned}
\end{equation} 
\end{proposition}

\begin{proof}  
In order to conclude the statement, we show that
\begin{equation}
\label{eqtoprove-}
\begin{aligned}  
T_{K \cap C}(x) \cap F(x) & \neq \emptyset ~~~ \forall x \in U(x_o) \cap \partial (K \cap C)~\mbox{and} \\ & ~~ \forall x_o \in (K \cap \partial C) \backslash D.
\end{aligned}
\end{equation}
Indeed, this latter fact, using \ifitsdraft Proposition \ref{thm2} \else 
\cite[Proposition 8]{draftautomatica}\fi, implies the existence of a nontrivial flow starting from each $x_o \in (K \cap \partial C) \backslash D$ which in turn allows to conclude the contractivity of the set $K$ using Definition \ref{def.cocopre-}. Now, in order to prove \eqref{eqtoprove-}, we distinguish three complementary situations. First, when $x \in  \partial (K \cap C) \cap \mbox{int}(C) $, in this case, using the forward pre-invariance of the set $K$, see Lemma \ref{lemmmm1}, and the fact that there exists always a nontrivial flow starting from $x \in \mbox{int}(C)$ under the standing assumptions, we conclude the existence of a nontrivial flow starting from $x$ and remaining in $K$. Hence, using \ifitsdraft Proposition \ref{thm1} \else 
\cite[Proposition 7]{draftautomatica}\fi, \eqref{eqtoprove-} follows when $x \in \partial (K \cap C) \cap \mbox{int}(C)$. Second, when $x \in  \partial (K \cap C) \cap \partial C  \cap \mbox{int}(K_e) = \partial C \cap \mbox{int}(K_e) $. In this case, \eqref{eqtoprove-} follows from \eqref{eqexists} since in this case 
$T_C(x) = T_{K \cap C} (x)$. Third, when $x \in  \partial (K \cap C) \cap \partial C \cap \partial K_e = \partial C \cap \partial K_e$, in this case under \eqref{eqexists}, we conclude the existence of $\eta \in F(x)$ such that $\eta \in T_C(x)$. Furthermore, using \eqref{eq.mink1pre} and the convexity of $K$, we conclude that $\eta \in D_{\mbox{int}(K)}(x)$, see the proof of \ifitsdraft Lemma \ref{lemmmm1} \else \cite[Lemma 1]{draftautomatica}\fi. Hence, using \ifitsdraft Lemma \ref{lem6} \else \cite[Lemma 5]{draftautomatica}\fi, it follows that $ \eta \in T_{K \cap C}(x) $. 
\ifitsdraft 
\else
\hfill $\square$
\fi
\end{proof}

\begin{example} \label{expCsets-}
We propose to build upon the pre-contractivity conclusions in Example \ref{expCsets} in order to conclude contractivity using Proposition \ref{lemconta}. To do so, it is enough to show that the set $(K \cap \partial C) \backslash D$ satisfies \eqref{eqexists}. Indeed, the set $(K \cap \partial C) \backslash D$ reduces to the singleton $\left\{ [-1~~-1]^\top \right\}$ and one can take the neighborhood $ U(x) \cap K \cap \partial C = \left\{ x \in \mathbb{R}^2 : x_2 = -1,~ x_1 \in [-1,-1/2)  \right\} $ on which $F(x) = [- [1,2] x_1 - 1/2 ~~ 1 - (1/2) x_1 ]^\top \in T_C(x)$  since its first component is $F_1(x) = 1 - (1/2) x_1 > 0$, hence, \eqref{eqexists} follows.
\end{example}

\subsection{Auxiliary Results}

\begin{lemma} \label{lem1}
Let $x$ be a nontrivial solution to $\mathcal{H}$ such that $[t'_1,t'_2] \times \left\{ j \right\} \subset \dom x$ for some $j \in \mathbb{N}$. Then, $ \dot{x}(t,j) \in T_C(x(t,j)) $ for almost all $t \in [t'_1, t'_2]$.
\end{lemma}
\begin{proof}
Let $t \in (t'_1, t'_2)$ such that $\dot{x}(t,j)$ exists thus $\dot{x}(t,j) \in F(x(t,j))$. Moreover, let a sequence $\left\{t_n\right\}_{n \in \mathbb{N}} \subset (0, t'_2 - t) $ such that $t_n \rightarrow 0$. That is, for $v_n(t) := (x(t_n,j) - x(t,j)) / t_n$, we have $\lim_n v_n(t) = \dot{x}(t,j)$ and at the same time $ x(t,j) + t_n v_n(t) = x(t_n,j) \in C $. Hence, using \eqref{eq.conti}, we conclude that 
$\dot{x}(t,j) \in T_C(x(t,j)) $.
\ifitsdraft 
\else
\hfill $\square$
\fi
\end{proof}

The statement of the following Lemma can also be found in \cite[Proposition 4.3.7]{aubin2009set}, however, the proof we propose in this paper is original and relatively simpler.

\begin{lemma} \label{lem3}
Consider the set $M \subset \mathbb{R}^n$ and the closed $K \subset M$, and a multiple barrier function candidate $B: \mathbb{R}^n \rightarrow \mathbb{R}^m$ such that $ K := \left\{ x \in M :~ B(x) \leq 0 \right\} $. Let $x \in \partial K \cap \mbox{int}(M)$ such that $B_i(x) = 0$ if and only if $i \in I_x \subset \left\{1,...,m\right\}$. Assume further that
\begin{enumerate}[label={(\roman*)},leftmargin=*]
\item \label{item:ij)}  there exists a neighborhood $U(x)$ such that $B_i$ is $\mathcal{C}^1(U(x))$ for all $i \in I_x $, 
\item \label{item:iij)} there exists $v \in \mathbb{R}^n$ such that $\nabla B_i(x)^\top v < 0$ for all $i \in I_x $. 
\end{enumerate}
Then
$$ T_K(x) = \left\{ w \in \mathbb{R}^n : \nabla B_i(x)^\top w \leq 0, ~~\forall i \in I_x \right\}. $$
\end{lemma}
\begin{proof} 
In the first step, we propose to show that for $x \in \partial K \cap \mbox{int}(M)$ and $w \in \mathbb{R}^n$, if $\langle \nabla B_i(x), w \rangle \leq 0$ for all $i \in I_x$ then $w \in T_K(x)$. To this end, we introduce the convex combination $ w_\beta := \beta v + (1-\beta) w$ with $v$ introduced in \ref{item:iij)}. That is, it is easy to see that, for any $\beta \in (0, 1]$, \ref{item:iij)} remains satisfied when we replace $v$ therein is replaced by $w_\beta$. Furthermore, since $x \in \mbox{int}(M)$ we claim that $w_\beta \in T_{K}(x) $ for all $\beta \in (0,1]$. Finally, using the closeness of the contingent cone $T_K(x)$ we conclude that $w = \lim_{\beta \rightarrow 0} w_\beta \in T_K(x) $. Now, in order to prove the claim, we use the Taylor expansion of $B_i(x + h w_\beta)$, for all $i \in I_x$, given by
\begin{align*} 
B_i(x + h w_\beta)  = \langle \nabla B_i(x), w_\beta \rangle h + h \epsilon_i(h),  
\end{align*}
where $\epsilon_i(h) \rightarrow 0$ as $h \rightarrow 0$. The latter expansion implies the existence of 
$\alpha > 0$ sufficiently small such that $B_i(x + h w_\beta) \leq 0$ for all $h \in (0, \alpha]$ and for all $i \in I_x$. Moreover, when $i \notin I_x$, we already have $B(x) < 0$, hence under the continuity of 
$B$ and for $\alpha$ sufficiently small, we conclude that $B(x + h w_\beta) \leq 0$ for all $h \in (0, \alpha]$. Furthermore, since $x \in \partial K \cap \mbox{int}(M)$, we conclude the existence of $\alpha_1 > 0$ such that, for each $h \in (0, \alpha_1]$, $ x + h w_\beta \in M$. Thus, for each $h \in (0,\min \left\{ \alpha, \alpha_1 \right\}]$, $x+h w_\beta \in K$, hence, $v \in T_{K}(x)$ using \eqref{eq.conti}.

In the second step, we show that if $w \in T_K(x)$ then $\nabla B_i(x)^\top v \leq 0$ for all $i \in I_x$. Indeed, having $K = M \cap K_e$, implies that $w \in T_{K_e}(x)$. Furthermore, having 
$$K_{e} = \cap^{m}_{i=1} K_{ei}, ~~~ K_{ei} = \left\{ x \in \mathbb{R}^n : B_i(x) \leq 0 \right\}, $$  
we conclude that $w \in \cap^{m}_{i=1} T_{K_{ei}}(x)$. Moreover, since $T_{K_{ei}}(x) = \mathbb{R}^n$ for all $i \notin I_x$, it follows that $w \in \cap_{i \in I_x} T_{K_{ei}}(x)$. Finally, we show that, under \ref{item:ij)}-\ref{item:iij)}, $ \nabla B_i (x)^\top v \leq 0 $ for all $i \in I_x$  using contradiction. That is, for $i \in I_x$, assume that $v \in T_{K_{ei}}(x)$ and at the same time $\langle \nabla B_i(x), v \rangle > 0$. Having $v \in T_{K_{ei}}(x)$ is equivalent to the existence of a sequences $h_i \rightarrow 0^+$ and $v_i \rightarrow v$ such that $x + h_i v_i \in K_{ei}$. Furthermore, using the Taylor expansion of $B_i(x + h_i v_i)$ and using the fact that $B_i(x) = 0$ for all $x \in \partial K_{ei}$, we obtain  
\begin{align*} 
B_i(x + h_i v_i)  = & \langle \nabla B_i(x), v \rangle h_i + \langle \nabla B_i(x), (v_i-v) \rangle h_i + \\ & h_i \epsilon_i(h_i) = 0,
\end{align*}
where $\epsilon(h) \rightarrow 0$ as $h \rightarrow 0$. However, since $ \langle \nabla B_i(x), v \rangle > 0$ we conclude the existence of $j \in \mathbb{N}$ such that, for each $i \geq j$, $B_i(x + h_i v_i) > 0$. Hence, $ x + h_i v_i \in \mathbb{R}^n \backslash K_{ei}$, which yields to a contradiction.
\ifitsdraft 
\else
\hfill $\square$
\fi
\end{proof} 
The following corollary is a particular case of Lemma \ref{lem3} where the barrier function candidate is scalar. Furthermore, it provides a characterization of $T_K$ only at the elements of $\partial K$ that are also in the interior of the set $M$. Similar statement can be found in \cite[Lemma 2.20]{prajna2005optimization}.
\begin{corollary} \label{lem2}
Consider closed sets $M \subset \mathbb{R}^n$ and $K \subset M$, let $B : \mathbb{R}^n \rightarrow \mathbb{R}$ a barrier function candidate defining the set $K \subset M$, 
$ K := \left\{ x \in M :~ B(x) \leq 0 \right\} $, and let $x \in \partial K \cap \mbox{int}(M)$ such that 
\begin{enumerate}[label={(\roman*)},leftmargin=*]
\item \label{item:ijc)}  $B$ is $\mathcal{C}^1$ in an open neighborhood around $x$ denoted $U(x)$,
\item \label{item:iijc)} $\nabla B(x) \neq 0$. 
\end{enumerate}
Then
\begin{align} \label{eq.conee} 
T_K(x) = \left\{ v \in \mathbb{R}^n : \nabla B(x)^\top v \leq 0 \right\}. 
\end{align}
\end{corollary}

\begin{lemma} \label{lem2bis}
Consider a closed set $M \subset \mathbb{R}^n$ and let $B : \mathbb{R}^n \rightarrow \mathbb{R}^m $ be a multiple barrier function candidate defining the set $K \subset M$ as $ K := 
\left\{ x \in M :~ B(x) \leq 0 \right\}$. Suppose $x \in \partial K \cap \mbox{int}(M)$ is such that
\begin{enumerate} [label={(\roman*)},leftmargin=*]
\item  $B$ is $\mathcal{C}^1$ in an open neighborhood around $x$ denoted $U(x)$,\\  
\item $\nabla B_i(x) \neq 0 ~\forall i \in I_x$, 
\begin{align*} 
I_x := \left\{ i \in \left\{1,..,m\right\} :  B_i(x) = 0 \right\}.
\end{align*}
\end{enumerate}
Then
$$ D_{\mbox{int}(K)}(x) = \left\{ v \in \mathbb{R}^n : \nabla B_i(x)^\top v < 0~~\forall i \in I_x \right\}. $$ 
\end{lemma}
\begin{proof}
We start noticing that since $x \in \mbox{int}(M)$, then, $D_{\mbox{int}(K)}(x) = D_{\mbox{int}(K_e \cap M)}(x)$. Moreover, using \eqref{eq.cone3} and Lemma \ref{lem6}, we obtain 
$$
D_{\mbox{int}(K_e \cap M)}(x) = D_{\mbox{int}(K_e)}(x) \cap D_{\mbox{int}(M)}(x) = D_{\mbox{int}(K_e)}(x).
$$ 
Furthermore, since $K_e = \cap^m_{i=1} K_{ei}$, it follows by repeating the previous reasoning that 
$$ D_{\mbox{int}(K)}(x) = D_{\mbox{int}(\cap^m_{i=1} K_{ei})}(x) = \cap^m_{i=1} D_{\mbox{int}(K_{ei})}(x). $$
Moreover, since $D_{\mbox{int}(K_{ei})}(x) = \mathbb{R}^n$ for all $i \notin I_x$, we obtain:
$$ D_{\mbox{int}(K)}(x) = \cap_{i \in I_x} D_{\mbox{int}(K_{ei})}(x). $$
Hence, the lemma is proved if we show that, for all $i \in I_x$,
\begin{align} \label{eq-eq}
D_{\mbox{int}(K_{ei})}(x) = & \left\{ v \in \mathbb{R}^n : \nabla B_i(x)^\top v < 0   \right\}. 
\end{align}
Indeed, the proof of \eqref{eq-eq} is obtained by combining Corollary \ref{lem2} with \eqref{eq.cone3}. Indeed, let us introduce the closed set $\mathbb{R}^n \backslash \mbox{int}(K_{ei})$. It is easy to see that 
$$ \mathbb{R}^n \backslash \mbox{int}(K_{ei}) = \left\{ y \in \mathbb{R}^n : \bar{B}(y) := -B(y) \leq 0 \right\}.  $$ 
Hence, using Corollary \ref{lem2} under the fact that $ x \in \partial (\mathbb{R}^n \backslash \mbox{int}(K_{ei})) $, we conclude that 
$$ T_{\mathbb{R}^n \backslash \mbox{int}(K_{ei})}(x) = \left\{ v \in \mathbb{R}^n : \nabla B(x) v \geq 0 \right\}. $$
Finally, applying \eqref{eq.cone3}, the statement follows.
\ifitsdraft 
\else
\hfill $\square$
\fi
\end{proof}

\begin{lemma} \label{lem6}
Let $v \in \mathbb{R}^n$ satisfying $v \in D_{\mbox{int}(K_1)}(x) \cap T_{K_2}(x) $ for $K_1$ and $K_2$ closed subsets of $\mathbb{R}^n$ and $x \in K_1 \cap K_2$. Then, $v \in T_{K_1 \cap K_2}(x)$. Moreover, if $v \in D_{\mbox{int}(K_2)}(x)$, then $v \in D_{\mbox{int}(K_1 \cap K_2)}(x)$.
\end{lemma}
\begin{proof}
The proof follows using the characterization of the contingent cone in \eqref{eq.conti}. Indeed, $v \in T_{K_2}(x)$ if and only if there exists sequences $h_i \rightarrow 0^+$ and $v_i \rightarrow v$ such that $x + h_i v_i \in K_2$. On the other hand, using \eqref{eq.cone2}, we conclude the existence of 
$\alpha > 0$ and $\epsilon > 0$, such that $x + (0, \alpha](v+\epsilon \mbox{int}(\mathbb{B}) ) \subset K_1$. Which implies the existence of $j \in \mathbb{N}$ such that for each $i \geq j$,  $x + h_i v_i \in x + (0, \alpha](v+\epsilon \mbox{int}(\mathbb{B}) ) \subset K_1 $, hence, $x + h_i v_i \in K_1 \cap K_2$ for all $i \geq j$, which implies that $ v \in T_{K_1 \cap K_2}(x)$. Furthermore, $v \in D_{\mbox{int}(K_2)}(x)$ implies the existence $\alpha_1>0$ and $\epsilon_1>0$ such that $x + (0, \alpha_1](v+\epsilon_1 \mbox{int}(\mathbb{B}) ) \subset K_2$. Hence, $ x + (0, \min\left\{\alpha, \alpha_1\right\}](v+ \min\left\{\epsilon, \epsilon_1\right\} \mbox{int}(\mathbb{B}) ) \subset K_2 \cap K_1 $ which implies that $v \in D_{\mbox{int}(K_1 \cap K_2)}(x)$.
\ifitsdraft 
\else
\hfill $\square$
\fi
\end{proof} 

\subsection{Background from the literature}

We start this section by recalling from \cite[Definition 2.6]{goebel2012hybrid} the concept of solutions to a hybrid system $\mathcal{H}$.
\begin{definition} [solution to $\HS$] \label{defsol}
\label{solution definition}
\index{solution to a hybrid system}
A function $x : \dom x \to \reals^n$ defined on a hybrid time domain  $\dom x$ and such that, for each $j \in \nats$, $t \mapsto x(t,j)$ is absolutely continuous is a {\em solution} to $\HS$ if
\begin{itemize}               
\item[(S0)] $x(0,0) \in \mbox{cl}(C) \cup D$;
\item[(S1)] for all $j \in \nats$ such that 
$I^j:=\defset{t}{(t,j)\in\dom x}$
has nonempty interior
\begin{equation}
\label{S1}
\begin{array}{lcl}               
x(t,j)\in C & \mbox{for all} & t\in \mbox{int}(I^j), \cr
\dot{x}(t,j)\in F(x(t,j)) & ~~\mbox{for a.a.} & t\in I^j; \cr
\end{array} 
\end{equation}                          
\item[(S2)] for all $(t,j)\in \dom x$ such that $(t,j+1)\in \dom x$,\begin{equation}                 
\label{S2}        
\begin{array}{l}         
x(t,j)\in D, \qquad
x(t,j+1)\in G(x(t,j)). \cr
\end{array}
\end{equation}                 
\end{itemize}     
\end{definition}
Next, we recall the regularity and the continuity notions used in this paper for a set-valued map $F: O \rightrightarrows \mathbb{R}^n$ with $O \subset \mathbb{R}^n$.
\begin{definition} [Continuity notions] \label{deflusc}
$F$ is \textit{upper semicontinuous} at $x \in O$ if for any neighborhood $U(F(x))$ there exists $\epsilon > 0$ such that, for any $y \in x + \epsilon \mbox{int}(\mathbb{B}) $, $F(y) \cap U(F(x))$. It is \textit{upper hemicontinuous} at $x \in O$ if, for any $p \in \mathbb{R}^n$, the single-valued map $y \mapsto \delta(F(y),p) := \sup_{z \in F(y)} \langle p, z \rangle \in (-\infty, +\infty]$ is upper semicontinuous at $x$. Furthermore, it is said to be upper semicontinuous or upper hemicontinuous, respectively, if it is so for all $x \in O$.
\end{definition}

\begin{definition} [Local boundedness]
A set-valued map $F: O \rightrightarrows \mathbb{R}^n$ is said to be \textit{locally bounded} if for any $x \in O$ there exists $U(x)$ and $K > 0$ such that 
\begin{align} \label{eq:lb}
|\zeta| \leq K ~~~~~~ \forall \zeta \in F(y)~~~\forall y \in U(x).  
\end{align}
\end{definition}

\begin{definition} [Lipschitz regularity] \label{deflip}
A set-valued map $F : O \rightrightarrows \mathbb{R}^m$ is said to be locally Lipschitz if for any compact set 
$V \subset O$ there exists $k>0$ such that 
\begin{align} \label{eq.lipset}
F(y) \subset & ~F(x) + k |x-y| \mbox{int}(\mathbb{B}) ~~~~~~\forall (x,y) \in V \times V.  
\end{align}
\end{definition}

On the other hand, we recall from the existing literature of differential inclusions some useful results that will play an important role in proving our statements. That is, we recall from \cite[Propositions 3.4.1 and 3.4.2]{Aubin:1991:VT:120830} the following statements that allow the extension of the well-know Nagumo's invariance theorem under differential inclusions. Those results will play an important role precisely when providing sufficient conditions for pre-contractivity and forward invariance using barrier functions.

\begin{proposition} \cite[Proposition 3.4.1]{Aubin:1991:VT:120830} \label{thm1}
Let us assume that the set-valued map $F: O \rightrightarrows \mathbb{R}^n$, with $O \subset \mathbb{R}^n$, satisfies 
\begin{enumerate}[label={(\roman*)},leftmargin=*]
\item \label{item:difinc}
 $F$ is upper hemicontinuous on $O$, 
\item \label{item:difinc1} $F(x)$ is convex and compact for all $x \in C \subset O$. 
\end{enumerate}
Consider a solution $x$ to the differential inclusion $\dot{x} \in F(x)$ starting at $x_o$ and satisfying 
$$ \forall T > 0, ~ \exists t \in (0, T] : x(t) \in K \subset C. $$  
Then,  $F(x_o) \cap T_K(x_o) \neq \emptyset$.
\end{proposition}

\begin{proposition} \cite[Proposition 3.4.2]{Aubin:1991:VT:120830} \label{thm2}
Let us assume that the set-valued map $F: O \rightrightarrows \mathbb{R}^n$ with $O \subset \mathbb{R}^n$ satisfies 
\begin{enumerate}[label={(\roman*)},leftmargin=*]
\item \label{item:1difinc} $F$ is upper semicontinuous on $O$, 
\item \label{item:2difinc} $F(x)$ is convex and compact for all $x \in C \subset O$. 
\end{enumerate}
Let the set $K \subset C$ be locally compact and let $K_o \subset K$ be a compact neighborhood of $x_o \in K$ such that
$$ \forall x \in K_o, ~~F(x) \cap T_K(x) \neq \emptyset.  $$  
Then, there exists $T>0$ and a solution to $\dot{x} \in F(x)$ starting at $x_o$ such that $x([0,T)) \subset K$.
\end{proposition}

\begin{remark} \label{rem-iinf}
it is useful to notice that closed subsets of finite dimensional spaces are locally compact.
\end{remark}

Furthermore, we recall the following result that will play an important role when proposing sufficient conditions for pre-contractivity using using barrier functions. 
\begin{theorem} \cite[Theorem 4.3.4]{Aubin:1991:VT:120830} \label{thm3}
Consider a nontrivial upper semicontinuous set-valued map $F: O \rightrightarrows \mathbb{R}^n$ such that $F(x)$ is convex and compact for any $x \in C \subset O$. Let $K \subset C$ be closed with nonempty interior and $x_o \in \partial K$. If $F(x_o) \subset D_{\mbox{int}(K)}(x_o)$, then, for each solution $x$ starting from $x_o$,
$$ \exists T > 0 : x((0,T]) \subset \mbox{int}(K). $$      
\end{theorem} 

\begin{remark} \label{remnose}
We stress that our standing assumptions on the flow map $F$ are equivalent to the regularities required in Propositions \ref{thm1}-\ref{thm2} and Theorem \ref{thm3}. Indeed, outer semicontinuous and locally bounded set-valued maps are upper semicontinuous with compact images \cite[Theorem 5.19]{rockafellar2009variational}, the converse is also true using \cite[Lemma 5.15]{goebel2012hybrid} and the fact that upper semicontinuous set-valued maps with compact images are locally bounded. Furthermore, upper semicontinuous set-valued maps are also upper hemicontinuous provided that the images are closed, see \cite[Corollary 2.4.1]{Aubin:1991:VT:120830}. The converse is true provided that the images are closed and convex, 
see \cite[Remark. Page 67]{Aubin:1991:VT:120830}.  
\end{remark} 

When a scalar function $B : \mathbb{R}^n \rightarrow \mathbb{R}$ is locally Lipschitz, its \textit{generalized gradient}, denoted by $\partial_C B$, constitutes a useful tool to understand its behavior along the system's solutions. We recall the following definition which is valid due to the equivalence established in \cite[Theorem 8.1]{clarke2008nonsmooth}. 

\begin{definition} \label{defgen}
Let $B : \mathbb{R}^n \rightarrow \mathbb{R}$ be locally Lipschitz. Let $\Omega$ be any subset of zero measure in $\mathbb{R}^n$, and let $\Omega_B$ be the set of points in $\mathbb{R}^n$ at which $B$ fails to be differentiable. Then
\begin{align} \label{eq.gg}
\partial_C B(x) := \mbox{co} \left\{ \lim_{i \rightarrow \infty} \nabla B(x_i) : x_i \rightarrow x,~x_i \notin \Omega_B,~x_i \notin \Omega \right\}.
\end{align}
\end{definition} 

\fi

\ifitsdraft
\else

\section*{Biography}

\begin{wrapfigure}{l}{25mm} 
\includegraphics[width=1in,height=1.25in,clip,
keepaspectratio]{MM}
\end{wrapfigure} 
\par \textbf{Mohamed Maghenem} received his Control-Engineer degree from the Polytechnical School of Algiers, Algeria, in 2013, his M.S. and Ph.D. degrees in Control from the University of Paris-Saclay, France, in 2014 and 2017, respectively. He is currently a Postdoctoral Fellow at the Electrical and Computer Engineering Department at the University of California at Santa Cruz. His research interests include dynamical systems theory (stability, safety, reachability, robustness, and synchronization), control systems theory (adaptive, time-varying, linear, non-linear, hybrid, robust, etc.) with applications to power systems, mechanical systems, and cyber-physical systems.
\par

\begin{wrapfigure}{l}{25mm} 
\includegraphics[width=1in,height=1.25in,clip,keepaspectratio]{RS}
\end{wrapfigure} 
\par \textbf{Ricardo. G. Sanfelice} received the B.S. degree in Electronics Engineering from the Universidad de Mar del Plata, Buenos Aires, Argentina, in 2001, and the M.S. and Ph.D. degrees in Electrical and Computer Engineering from the University of California, Santa Barbara, CA, USA, in 2004 and 2007, respectively. In 2007 and 2008, he held postdoctoral positions at the Laboratory for Information and Decision Systems at the Massachusetts Institute of Technology and at the Centre Automatique et Systèmes at the École de Mines de Paris. In 2009, he joined the faculty of the Department of Aerospace and Mechanical Engineering at the University of Arizona, Tucson, AZ, USA, where he was an Assistant Professor. In 2014, he joined the University of California, Santa Cruz, CA, USA, where he is currently Professor in the Department of Electrical and Computer Engineering. Prof. Sanfelice is the recipient of the 2013 SIAM Control and Systems Theory Prize, the National Science Foundation CAREER award, the Air Force Young Investigator Research Award, the 2010 IEEE Control Systems Magazine Outstanding Paper Award, and the 2020 Test-of-Time Award from the Hybrid Systems: Computation and Control Conference. His research interests are in modeling, stability, robust control, observer design, and simulation of nonlinear and hybrid systems with applications to power systems, aerospace, and biology.
\par

\fi

\balance

\end{document}

%% file: Example5.pdf_tex
\begingroup%
  \makeatletter%
  \providecommand\color[2][]{%
    \errmessage{(Inkscape) Color is used for the text in Inkscape, but the package 'color.sty' is not loaded}%
    \renewcommand\color[2][]{}%
  }%
  \providecommand\transparent[1]{%
    \errmessage{(Inkscape) Transparency is used (non-zero) for the text in Inkscape, but the package 'transparent.sty' is not loaded}%
    \renewcommand\transparent[1]{}%
  }%
  \providecommand\rotatebox[2]{#2}%
  \newcommand*\fsize{\dimexpr\f@size pt\relax}%
  \newcommand*\lineheight[1]{\fontsize{\fsize}{#1\fsize}\selectfont}%
  \ifx\svgwidth\undefined%
    \setlength{\unitlength}{582.31521776bp}%
    \ifx\svgscale\undefined%
      \relax%
    \else%
      \setlength{\unitlength}{\unitlength * \real{\svgscale}}%
    \fi%
  \else%
    \setlength{\unitlength}{\svgwidth}%
  \fi%
  \global\let\svgwidth\undefined%
  \global\let\svgscale\undefined%
  \makeatother%
  \begin{picture}(1,0.38631573)%
    \lineheight{1}%
    \setlength\tabcolsep{0pt}%
    \put(0,0){\includegraphics[width=\unitlength,page=1]{Example5.pdf}}%
    \put(0.14322345,0.25880542){\color[rgb]{0,0,0}\makebox(0,0)[lt]{\lineheight{1.25}\smash{\begin{tabular}[t]{l}$F(x)$\end{tabular}}}}%
    \put(0.39226625,0.13368949){\color[rgb]{0,0,0}\makebox(0,0)[lt]{\lineheight{1.25}\smash{\begin{tabular}[t]{l}$x_1$\end{tabular}}}}%
    \put(0.16075313,0.35365135){\color[rgb]{0,0,0}\makebox(0,0)[lt]{\lineheight{1.25}\smash{\begin{tabular}[t]{l}$x_2$\end{tabular}}}}%
    \put(0,0){\includegraphics[width=\unitlength,page=2]{Example5.pdf}}%
    \put(0.22966304,0.06575023){\color[rgb]{0,0,1}\makebox(0,0)[lt]{\lineheight{1.25}\smash{\begin{tabular}[t]{l}$C$\end{tabular}}}}%
    \put(0.01633549,0.11629504){\color[rgb]{1,0,0}\makebox(0,0)[lt]{\lineheight{1.25}\smash{\begin{tabular}[t]{l}$K_e$\end{tabular}}}}%
    \put(0,0){\includegraphics[width=\unitlength,page=3]{Example5.pdf}}%
    \put(0.64960058,0.25751373){\color[rgb]{0,0,0}\makebox(0,0)[lt]{\lineheight{1.25}\smash{\begin{tabular}[t]{l}$F(x)$\end{tabular}}}}%
    \put(0.89864338,0.1323978){\color[rgb]{0,0,0}\makebox(0,0)[lt]{\lineheight{1.25}\smash{\begin{tabular}[t]{l}$x_1$\end{tabular}}}}%
    \put(0.66713026,0.35235966){\color[rgb]{0,0,0}\makebox(0,0)[lt]{\lineheight{1.25}\smash{\begin{tabular}[t]{l}$x_2$\end{tabular}}}}%
    \put(0,0){\includegraphics[width=\unitlength,page=4]{Example5.pdf}}%
    \put(0.72472724,0.05691659){\color[rgb]{0,0,1}\makebox(0,0)[lt]{\lineheight{1.25}\smash{\begin{tabular}[t]{l}$C_1$\end{tabular}}}}%
    \put(0.52271262,0.11500335){\color[rgb]{1,0,0}\makebox(0,0)[lt]{\lineheight{1.25}\smash{\begin{tabular}[t]{l}$K_e$\end{tabular}}}}%
  \end{picture}%
\endgroup%